\documentclass{amsart}
\usepackage{import}
\usepackage[latin9]{inputenc}
\usepackage{microtype}
\usepackage{fancyhdr}
\usepackage{cite} 
\usepackage{amsfonts}
\usepackage{amsmath}
\usepackage{amsthm}
\usepackage{amssymb}
\usepackage[all]{xy}
\usepackage{mathrsfs}
\usepackage{bm} 
\usepackage{stmaryrd}
\usepackage{color}   
\usepackage[hidelinks]{hyperref}
\hypersetup{
    colorlinks=false, 
    linktoc=all,     
    linkcolor=black,  
}
\usepackage{graphicx}
\usepackage{wrapfig}
\usepackage{float}
\usepackage{caption}
\usepackage{subcaption}
\usepackage{enumitem}
\usepackage{empheq} 
\usepackage{mathtools} 
\usepackage{adjustbox}
\usepackage{tikz,tikz-cd}
\usetikzlibrary{positioning,arrows}
\setlist{nosep} 
\usepackage{dutchcal} 

\usepackage{makecell} 

\usepackage[foot]{amsaddr} 

\usepackage{imakeidx} 
\makeindex 

\usepackage[splitrule]{footmisc} 

\newcommand{\quotient}[2]{{\raisebox{0.2em}{$#1$}
           \!\!\left/\raisebox{-0.2em}{\!$#2$}\right.}}

\newcommand{\unicov}{\widetilde}

\newcommand{\inject}{\hookrightarrow}

\usepackage[top=2.5cm,bottom=2cm, margin=2.5cm]{geometry}
\usepackage[font=small, width=.75\textwidth]{caption}
\usepackage{mwe}

\theoremstyle{definition}
\newtheorem{theorem}{Theorem}

\theoremstyle{definition}
\newtheorem{lemma}[theorem]{Lemma}

\theoremstyle{definition}
\newtheorem{question}[theorem]{Question}

\theoremstyle{definition}
\newtheorem{definition}[theorem]{Definition}

\theoremstyle{definition}
\newtheorem{proposition}[theorem]{Proposition}

\theoremstyle{definition}

\theoremstyle{definition}
\newtheorem{corollary}[theorem]{Corollary}

\theoremstyle{definition}
\newtheorem{remark}[theorem]{Remark}

\theoremstyle{definition}

\newcommand{\innerthmname}{}

\theoremstyle{definition}

\usepackage{titlesec}
\titleformat{\section}
  {\normalfont\fontsize{14}{15}\bfseries}{\thesection}{1em}{}

\titleformat{\subsection}
  {\normalfont\fontsize{12}{15}\bfseries}{\thesubsection}{1em}{}

\usepackage{xcolor}

\newcommand*{\colorboxed}{}
\def\colorboxed#1#{%
  \colorboxedAux{#1}%
}
\newcommand*{\colorboxedAux}[3]{%
  \begingroup
    \colorlet{cb@saved}{.}%
    \color#1{#2}%
    \boxed{%
      \color{cb@saved}%
      #3%
    }%
  \endgroup
}

\title{\vspace{-3cm}Cusp Shape and Fractional Dehn Twists of\\ Fibred Hyperbolic 3-Manifolds}
\author{Misha Schmalian}
\address{\hskip -0.35cmMathematical Institute, University of Oxford, OX2 6GG, UK, https://sites.google.com/view/mishaschmalian}
\begin{document}
\vspace*{3cm}
\maketitle
\date{}
\vspace{-1cm}
\begin{abstract}
Given a fibred hyperbolic 3-manifold with boundary, we
coarsely relate the Euclidean geometry of its cusps to the classical
fractional Dehn twist coefficient of its monodromy. This result fits
into the broader programme of coarsely describing the geometry of a
hyperbolic 3-manifold via combinatorial data. We are thus able to study
the hyperbolic geometry of certain fibred 3-manifolds under Dehn
filling. For example, we find coarse volume estimates for sufficiently
twisted braid closures in terms of their braid words. We also prove that
for any open book decomposition of a fixed manifold (that is not a lens
space or solid torus) with fibre of fixed Euler characteristic the
fractional Dehn twist coefficient in some boundary component is uniformly bounded.
Finally, we obtain applications to contact topology. We give a geometric
criterion on the binding of an open book decomposition for the
corresponding contact structure in a hyperbolic 3-manifold to be tight.
\end{abstract}

\section{Introduction}\label{sec:1}
\noindent
Fix a finite type oriented surface $S$ with boundary. Consider an orientation-preserving diffeomorphism $\varphi:S\to S$ pointwise fixing the boundary $\partial S$ which is isotopic to a pseudo-Anosov diffeomorphism. The mapping torus of $\varphi$ is the following oriented quotient space:
$$M_\varphi := \,\quotient{\text{Int}(S) \times \mathbb R}{(x, t+1) \sim (\varphi(x), t)}.$$
Mapping tori form a large class of 3-manifolds, as evidenced, for example, by the virtual fibring theorem \cite{Ago13}. By \cite{Mos68, Thu86}, the mapping torus $M_\varphi$ admits a unique complete, finite-volume hyperbolic geometry. As is often the case, the existence proof of a hyperbolic geometry is non-constructive. However, as we discuss in section \ref{section:3}, substantial previous work coarsely relates geometric properties of $M_\varphi$ to combinatorial properties of $\varphi$. In this article, we study the cusp geometry of $M_\varphi$. \\

\noindent
Given a boundary component $R\subset \partial S$, the corresponding boundary component $T=R\times S^1\subset \partial M_\varphi$ has a neighbourhood $C_R=C$ called the maximal cusp, which we define in section \ref{section:3}. The hyperbolic metric on $M_\varphi$ restricts to a Euclidean metric on $\partial C$. This Euclidean geometry has been extensively studied and is closely related to the behaviour of the manifold under Dehn filling, as is surveyed in \cite{Lac19}. In this article we consider $S$ with boundary and punctures and let $\text{Diff}^+(S, \partial S)$ denote the set of diffeomorphisms of $S$ pointwise fixing the boundary $\partial S$, but not necessarily fixing the punctures. 

\begin{definition}\label{definition:1} In our setup, there is a standard choice of basis for $\pi_1(\partial C)$ consisting of $\gamma_{lat,R}=\gamma_{lat}:=\,\quotient{\{q\}\times \mathbb R}{\sim}$, for $q\in R$ and $\gamma_{\partial,R}=\gamma_\partial:=\,\quotient{R\times \{0\}}{\sim}$. The Euclidean torus $\partial C$ is orientation-preservingly isometric to $\mathbb R^2/\langle(a,0), (b,c)\rangle$ with $a,c>0$ and $(a,0), (b,c)$ corresponding to $\gamma_{\partial}$ and $\gamma_{lat}$ respectively. Consult figure \ref{figure:1}. We thus define the following real-valued quantities, which fully determine the
cusp geometry:
\begin{itemize}
\item[-] $l(\gamma_{\partial}):=a$; 
\item[-] $\text{height}(\partial C):=c=\text{Area}(\partial C)/l(\gamma_{\partial})$; 
\item[-] the {\it cusp-skew} $\mathcal{sk}(\varphi, R)=\mathcal{sk}(\varphi):=b/a$.
\end{itemize}
\end{definition}

\noindent
Informally, the cusp-skew is a measure of how ``slanted'' the induced Euclidean metric on $\partial C$ is. The quantities $l(\gamma_{\partial})$ and $\text{height}(\partial C)$ are well-understood by \cite{FS14}. Our main result is the following:

\begin{theorem}\label{theorem:2}
Consider an oriented, compact surface $S$, $\varphi\in \text{Diff}^+(S, \partial S)$ isotopic to a pseudo-Anosov map, a
boundary component $R \subset \partial S$, and the corresponding boundary
component $T=R\times S^1$ of the mapping torus $M_\varphi$. Then
the cusp-skew $\mathcal{sk}(\varphi, R)$ of the maximal cusp boundary $\partial C$
corresponding to $T$ satisfies:

\vskip-0.2cm
$$|\mathcal{fD}(\varphi, R)-\mathcal{sk}(\varphi, R)|\leq 6|\chi(S)|\cdot d_{\mathcal{A}(S, R)}(\varphi)+3.$$
\vskip 0.1cm

\noindent
Here $\mathcal{fD}(\varphi, R)$ denotes the fractional Dehn twist coefficient, a
well-studied real-valued quantity, which we define in section \ref{section:4}; $\chi(S)$ denotes the Euler characteristic; $\mathcal{A}(S,R)$
denotes the {\it restricted arc graph} which is the induced subgraph of
the (regular) arc graph consisting of arcs with at least one endpoint on
$R$; and finally $d_{\mathcal{A}(S,R)}(\varphi):=\inf_{v\in \mathcal{A}(S,R)}d(v, \varphi(v))$ is the translation distance in $\mathcal{A}(S, R)$.
\end{theorem}

{\centering
\includegraphics[width=8.5cm]{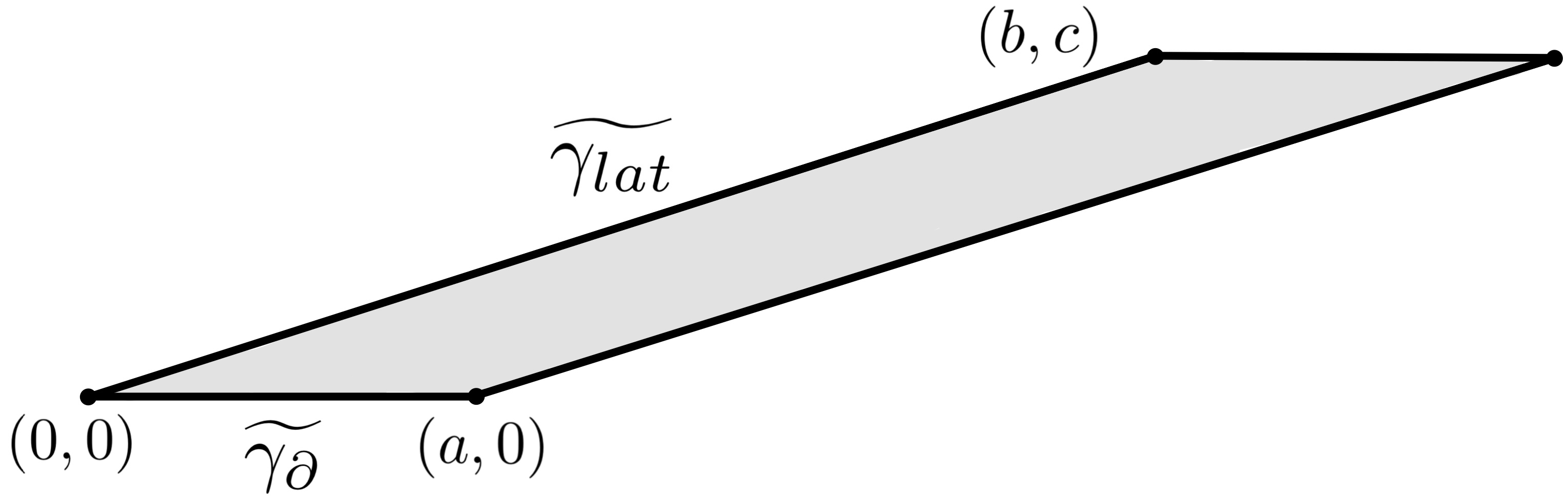}
\captionof{figure}{Lifts $\unicov{\gamma_{lat}}, \unicov{\gamma_\partial}$ in $\mathbb R^2$ for $\gamma_{lat}, \gamma_{\partial}$ as described in definition \ref{definition:1}.}
\label{figure:1}}
\vskip 0.3cm

\noindent
This completes the coarse combinatorial picture for the cusp geometry of $M_\varphi$. The article continues with several applications of this result.\\

\noindent
First of all, we observe that if $\mathcal{fD}(\varphi)$ is large, then we have good geometric behaviour under Dehn filling as evidenced by the following result:

\begin{theorem}\label{theorem:3} Let $S$ be a compact, oriented surface, $\varphi\in \text{Diff}^+(S, \partial S)$, $J>1, 0<\varepsilon\leq \log(3)$ and let $M$ be the Dehn filling of the mapping torus $M_\varphi$ along some subset of the curves $\{\gamma_{lat,R} |R\subset \partial S\}$. Suppose that:
\begin{itemize}
\item[-] For every boundary component $R\subset \partial S$ that we Dehn fill: $|\mathcal{fD}(\varphi, R)|>F(S, \varepsilon, J)$ for some constant $F=F(S, \varepsilon, J)$ given in remark \ref{remark:21}; 
\item[-] $\varphi$ is isotopic to a pseudo-Anosov map.
\end{itemize}
Then $M$ is hyperbolic and there are $J$-bi-Lipschitz inclusions $M^{\geq \varepsilon}\inject M_\varphi^{\geq \varepsilon/1.2}, M_\varphi^{\geq \varepsilon}\inject M^{\geq \varepsilon/1.2}$. (Here $M^{\geq\varepsilon}$ denotes the $\varepsilon$-thick part of $M$.)
\end{theorem}

\noindent A more concrete application is the geometry of braid closures. See section \ref{section:7} for a short introduction to braids. For a large class of hyperbolic braid closures we give a purely combinatorial formula for the volume, as follows:

\begin{theorem} \label{theorem:4} There exists $A=A(n)\in \mathbb R$ such that for a braid $\beta\in B_n$ with hyperbolic braid closure $\hat\beta$ and $|\lfloor \beta\rfloor_D|\geq 7\cdot (1+3216\cdot (n-1)^5)+4$ we have: 
$$\frac 1A\cdot \text{min}_{\rho\in \text{Conj}(\beta)}l_\mathcal{T}(\rho)\leq \text{vol}(S^3-\hat\beta)\leq A\cdot \text{min}_{\rho\in \text{Conj}(\beta)}l_\mathcal{T}(\rho).$$
Here $\text{Conj}(\beta)\subset B_n$ denotes the conjugacy class of $\beta$, $\lfloor \beta\rfloor_D$, defined in section \ref{section:7}, is the classical Dehornoy floor, and $l_\mathcal{T}(\rho)$, also defined in section \ref{section:7}, is the length of $\rho$ in terms of a certain generating set $\mathcal{T}$. 
\end{theorem}

\noindent 
Conversely, we may use hyperbolic geometry to study the fractional Dehn twist coefficient. We study all open book decompositions of a fixed 3-manifold and obtain the following result:

\begin{theorem}\label{theorem:5}
Given an oriented manifold $M$, that is not a lens space or a solid torus, there exists $N_{M, \chi}$ such that any abstract open book decomposition $(S, \varphi)$ of $M$ with $\chi(S)=\chi$ has $|\mathcal{fD}(\varphi, R)|\leq N_{M, \chi}+1$ for some component $R\subset \partial S$. 
\end{theorem}

\noindent
Here the assumption that $M$ is not a lens space or solid torus is necessary. This theorem, but with the assumption $\chi(S)=\chi$ removed and replaced by the assumption that $S$ has connected binding, was proved using Heegard Floer homology in \cite{HM18}. \\

\noindent
Finally, as we discuss in section \ref{section:9}, the fractional Dehn twist coefficient is related to contact topology. We may use this to relate the contact topology and hyperbolic geometry of 3-manifolds. As we mention in section \ref{section:4}, there is a dichotomy between so called tight and overtwisted contact structures and we will give geometric criteria for a contact structure to be tight.

\begin{theorem}\label{theorem:6} Let $\xi$ be a contact structure on an oriented closed hyperbolic 3-manifold $M$. Let $\gamma$ be the binding and $S$ the page of an open book decomposition given by the Giroux correspondence, as discussed in section \ref{section:4}. Suppose $\gamma$ is a geodesic with length $l$ and torsion $\tau$. (Informally the torsion $\tau$, defined in section \ref{section:9}, is a measure of how ``twisted'' $\gamma$ is.) There exist positive constants $L_{max} = L_{max}(S)$ and $tor_{min}=tor_{min}(S)$ depending only on $\chi(S)$, such that if $l^2+\tau^2\leq L_{max}, \tau/l\geq tor_{min}$, then $\xi$ is (positively) tight.
\end{theorem}

\noindent
In this article, unless stated otherwise, manifolds are 3-dimensional and oriented with, possibly empty, toroidal boundary. All manifolds are compact with some boundary components possibly removed. Manifolds are called hyperbolic if their interiors admit a complete, finite volume Riemannian metric of constant curvature $-1$.

\section{Acknowledgements} \label{section:2}
\noindent
The author would like to thank Marc Lackenby for all his input and
support, without which this paper would not have been possible. The author would also like to thank Jessica Purcell for her helpful comments regarding the angle-deformation theory used in this article.

\section{Background on Cusp Geometry} \label{section:3}
\noindent
It has become apparent that studying hyperbolic geometry is central to 3-manifold topology and many classes of 3-manifolds have been proven to admit complete, finite-volume hyperbolic metrics. However, often these proofs are non-constructive and questions remain about the geometric properties of hyperbolic 3-manifolds. Let us mention some previous results coarsely answering these question.

\begin{theorem}\label{theorem:7} \cite{Bro01b, Bro01a} Given a pseudo-Anosov diffeomorphism $\psi$ of a closed surface $S$ the mapping torus $M_\psi$ is known to be hyperbolic and has volume satisfying: 
$$d_\mathcal{P}(\psi)\leq K\text{vol}(M_\psi)\leq K^2 d_\mathcal{P}(\psi), $$
for some $K=K(S)>1$, with $\mathcal{P}$ the pants complex of $S$, and $d_\mathcal{P}(\psi):=\text{inf}_{v\in \mathcal{P}}d_{\mathcal{P}}(v, \psi(v)).$
\end{theorem}

\noindent
More generally, \cite{Min10} gives a construction of a piecewise Riemannian 3-manifold which is bi-Lipschitz, with constants depending on $S$, to $M_\psi$. \\

\noindent
In this article we will study the cusp geometry of $M_\varphi$. Substantial work on this topic has appeared in \cite{FS14} and we state their main result below. Let us first give an introduction to cusps of hyperbolic 3-manifolds.\\

\noindent
A boundary component $T\subset \partial M$ of a complete, finite-volume hyperbolic 3-manifold $M$ has a neighbourhood isometric to the quotient of $\{(x, y, z) | z\geq 1\}\subset \mathbb H^3$ in the upper-half space by a lattice $\Lambda \cong \mathbb Z^2$ acting by translation on the first two coordinates of $\mathbb H^3$. Such a neighbourhood, as sketched in figure \ref{figure:2}, is called a {\it cusp} and the union of all cusps at $T$ is called the {\it maximal cusp}. Given an element $\gamma\in \pi_1(T)$, let $l(\gamma)$ denote the length of any Euclidean geodesic representative of $\gamma$ in the boundary of the maximal cusp $C$ corresponding to $T$. As we mentioned above, the cusp geometry of a manifold is critical in understanding Dehn fillings and has been extensively studied. For example, the Dehn filling $M(\gamma)$ along a slope $\gamma$ with $l(\gamma)$ large enough remains hyperbolic and has similar hyperbolic geometry to $M$ \cite{Thu79}. Instead of $l(\gamma)$ in some contexts one considers the {\it normalised length} $\hat{L}(\gamma):=l(\gamma)/\sqrt{\text{Area}(\partial C)}$. We recall the following classical lower bound:  \\

{\centering
\includegraphics[width=6.3cm]{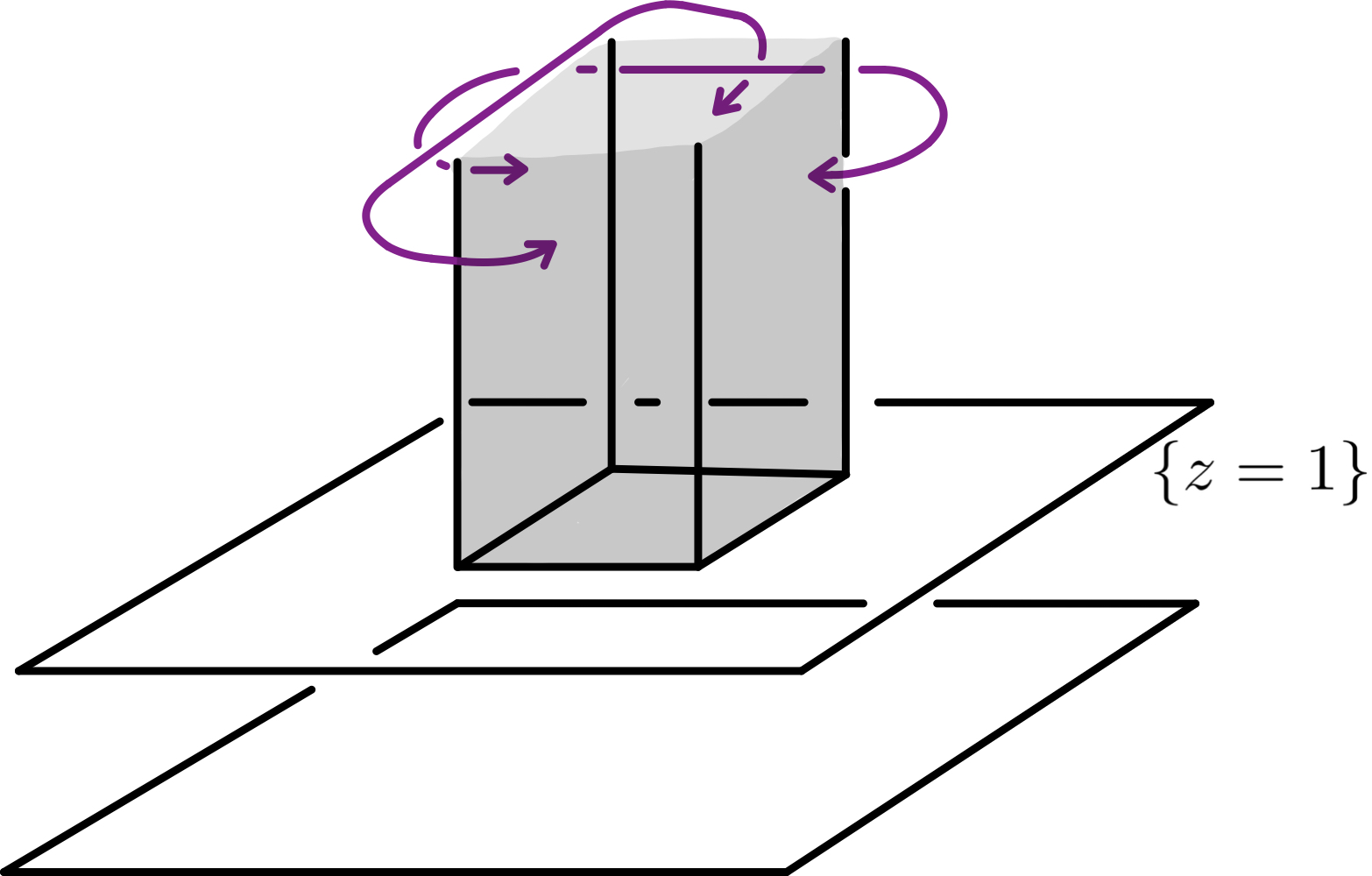}
\captionof{figure}{Fundamental domain of a cusp $C \subset M$ in the universal cover $\unicov{M}= \mathbb H^3$.}
\label{figure:2}}
\vskip 0.3cm

\begin{lemma} \label{lemma:8} In the above setup, for any non-trivial $\gamma\in \pi_1(T)$, we have $l(\gamma)\geq 1$. 
\end{lemma}

\noindent
{\it Proof:} For the following consult figure \ref{figure:3} left. We may take the universal cover of $M$ to be the upper half-plane model $\mathbb R^2\times \mathbb R_{>0}$ of $\mathbb H^3$ such that $C$ has a lift $\unicov{C_0}:=\{z\geq 1\}$. On $\mathbb R^2\times \mathbb R_{>0}$ we will consider both the hyperbolic and Euclidean metric. Since the cusp $C$ is maximal, there is another lift $\unicov{C_1}$ of $C$ touching $\unicov{C_0}$ at some point $\unicov{p}$. Take a lift $\unicov{\gamma}$ of the Euclidean geodesic representative of $\gamma$ in $C$ at $\unicov{p}$. There must be another lift $\unicov{C_2}$ of $C$ touching $\unicov{C_0}$ at the other endpoint $\unicov{q}$ of $\unicov{\gamma}$. Now $\unicov{C_1}, \unicov{C_2}\subset \mathbb R^2\times \mathbb R_{>0}$ are disjoint Euclidean balls of Euclidean diameter 1. Thus $\unicov{p}, \unicov{q}$ have Euclidean distance at least $1$. However, the metric on $\partial \unicov{C_0}$ induced by the hyperbolic metric and the Euclidean metric on $\mathbb R^2\times \mathbb R_{>0}$ agree. It follows that $l(\gamma)\geq 1$. $\square$\\

\noindent
Now suppose we have a proper incompressible embedded surface with boundary $S\subset M$ and an ideal triangulation $\tau$ of $S$. We may homotope the embedding $S\inject M$ to a map $f:S\to M$ such that $f(\tau)$ are geodesics and components of $S-\tau$ are mapped to totally geodesic surfaces. This $f$ is unique and is called the {\it pleating} of $S$ along $\tau$. The following upper bound is due to W. Thurston, see also \cite[Lemma 2.4]{FS14}:\\

{\centering
\includegraphics[width=14cm]{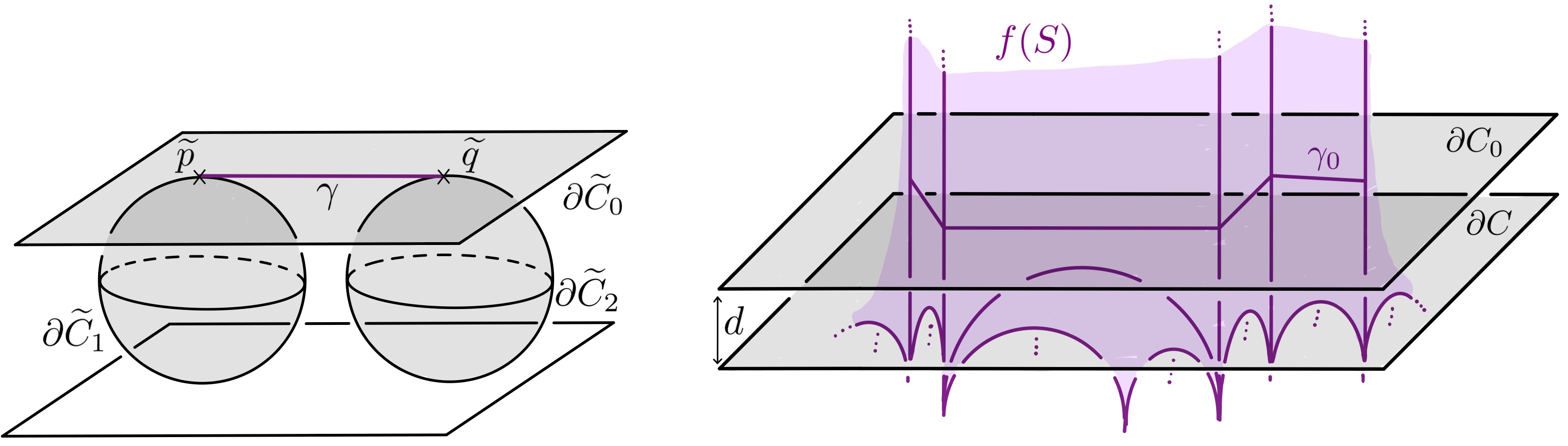}

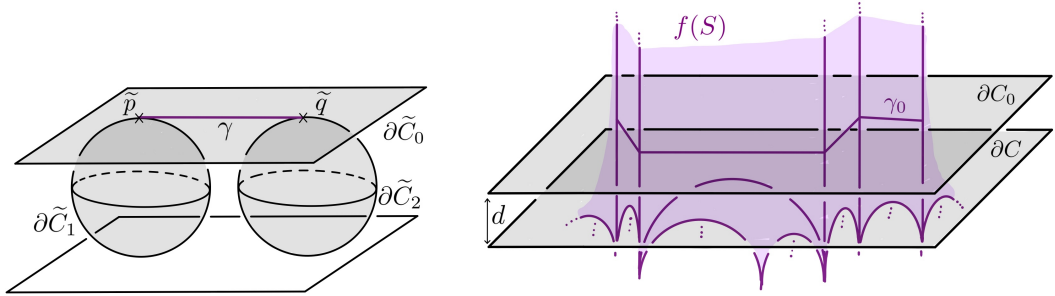
\captionof{figure}{Left: Lifts of cusp $C$ as described in proof of lemma \ref{lemma:8}; Right: Pleated surface $f(S)$ and cusps $C_0,C$ as described in proof of lemma \ref{lemma:9}.}
\label{figure:3}}
\vskip 0.3cm

\begin{lemma} \label{lemma:9} Consider the above setup and let $\gamma\in \pi_1(T)$ be the image of a loop about a puncture of $S$. Then $$l(\gamma)\leq 6|\chi(S)|.$$
\end{lemma}

\noindent
{\it Proof:} For the following consult figure \ref{figure:3} right. Recall, that $f$ is some pleating of $S$. Take a sufficiently small cusp $C_0\subset C$ of $M$, such that $f(S)\cap C_0$ consists of tips of ideal triangles and some component $\gamma_0\subset f(S)\cap \partial C_0$ is a closed curve homotopic to $\gamma$. Let $d$ denote the distance between $\partial C_0$ and $\partial C$ so that the Euclidean geometries on $\partial C_0$ and $\partial C$ are related by a dilation of factor $e^d$. In particular, $l(\gamma)\leq e^d\cdot \text{length}(\gamma_0)$. Equip $S$ with the hyperbolic metric pulled back along $f$, so that $f^{-1}(C_0)$ is a collection of cusps in $S$. Let $D_0\subset f^{-1}(C_0)$ be the cusp about the puncture defining $\gamma$. Observe that any geodesic in $M$ from $C_0$ to $C_0$ has length at least $2d$. Hence the cusp $D_0$ is contained in a larger embedded cusp $D\subset S$ with $\partial D_0$ and $\partial D$ a distance $d$ apart. Thus $\text{length}(\gamma_0) = \text{length}(\partial D_0) = e^{-d} \cdot \text{length}(\partial D)$. It is a standard fact that a 2-dimensional cusp has area equal to its boundary length. A horoball packing bound \cite{Bor78} states that $\text{area}(\emph{D})\leq \frac{6}{2\pi} \text{area}(S)$. Thus the Gauss-Bonet theorem gives the desired inequality $l(\gamma)\leq \text{length}(\partial D)=\text{area}(D)\leq \frac{6}{2\pi}\text{area}(S)=6|\chi(S)|$. $\square$\\

\noindent
Recall the setup from the introduction, that is: $S$ is an oriented surface with boundary and punctures, $\varphi\in \text{Diff}^+(S, \partial S)$ pointwise fixes the boundary but not necessarily the punctures, $R\subset \partial S$ is a boundary component, and $T\subset \partial M_\varphi$ is the corresponding boundary torus. Moreover, recall that with respect to the basis $\gamma_{lat}, \gamma_\partial$ the Euclidean geometry of the maximal cusp boundary $\partial C$ corresponding to $T$ is determined by the quantities $l(\gamma_\partial), \text{height}(\partial C)$, and $\mathcal{sk}(\varphi, R)$. It follows from lemmas \ref{lemma:8}, \ref{lemma:9}, that $1\leq l(\gamma_\partial)\leq 6|\chi(S)|$. Moreover, the following is known:

\begin{theorem}\label{theorem:10} \cite{FS14} In the above setup, we have: 
$$\frac{d_{\mathcal{A}(S, R)}(\varphi)}{536\chi(S)^4}<\text{height}(\partial C)\leq 3|\chi(S)|d_{\mathcal{A}(S, R)}(\varphi).$$
Therefore, to understand the cusp geometry of $M_\varphi$ it only remains to understand the cusp-skew $\mathcal{sk}(\varphi, R)$. 
\end{theorem}

\section{Background on Fractional Dehn Twist Coefficients} \label{section:4} 

\noindent
Consider the setup of $S$ an oriented surface with boundary components $R_1, \ldots, R_b$, $\chi(S)<0$, and $\varphi\in \text{Diff}^+(S, \partial S)$. Fix a hyperbolic metric on $S$ with geodesic boundary. The following definition is due to \cite{HKM07}:

\begin{definition}\label{definition:11} The {\it fractional Dehn twist coefficient} $\mathcal{fD}(\varphi, R_i)$ is defined as follows: 
\begin{itemize}
\item[-] Suppose $\varphi$ is freely isotopic to a pseudo-Anosov diffeomorphism $\hat\varphi$ preserving a geodesic lamination $L$. The component $C\subset S\setminus L$ containing $R_i$, as shown in figure \ref{figure:4}, is homeomorphic to an ideal $m$-gon for some $m\in \mathbb N$ with ideal vertices denoted $q_0, \ldots, q_{m-1}$ in clockwise order. Let $\alpha_k, k=0, \ldots, {m-1}$, as shown, be geodesic rays orthogonal to $R_i$ at an initial point $p_k\in R_i$ and tending towards $q_k$. Let $H:R_i\times [0,1]\to R_i$ be the restriction of the isotopy from $\varphi$ to $\hat\varphi$. Then exists $t\in \mathbb N, j\in \{0, \ldots, m-1\}$ such that $H(p_0, t), t\in [0,1]$ is isotopic, relative endpoints, to a clockwise path starting at $p_0$ going to $p_j$ and then winding $t$ full times around $R_i$. We define: 
$$\mathcal{fD}(\varphi, R_i):=t+\frac{j}{m};$$
\item[-] Suppose that $\varphi$ is freely isotopic to a periodic map. There exists $N$ and $t_1, \ldots, t_b\in \mathbb N$ such that $\varphi^N$ is isotopic, relative boundary, to a concatenation $T_{R_1}^{t_1}\circ \ldots \circ T_{R_b}^{t_b}$ of boundary Dehn twists. We define: $$\mathcal{fD}(\varphi, R_i):=\frac{t_i}{N}; $$
\item[-] Suppose $\varphi$ is reducible. Then there exists a subsurface $S'\subset S$ containing $R_i$ such that $\varphi(S')=S'$ and such that $\varphi|_{S'}$ is isotopic to a pseudo-Anosov or periodic map. We define $\mathcal{fD}(\varphi, R_i):=\mathcal{fD}(\varphi|_{S'}, R_i)$. 
\end{itemize}
\end{definition}

\noindent
Informally, $\mathcal{fD}(\varphi, R)$ indicates how much we ``twist'' $\varphi$ in $R$ to obtain its Nielsen-Thurston representative. Given $\varphi, \psi\in \text{Diff}^+(S, \partial S)$ and $n\in \mathbb Z$, the fractional Dehn twist coefficient satisfies the following properties: 
$$|\mathcal{fD}(\varphi\circ \psi, R)-\mathcal{fD}(\varphi, R)-\mathcal{fD}(\psi, R)|\leq 1, ~~~\mathcal{fD}(\varphi^n, R)=n\cdot \mathcal{fD}(\varphi, R),\text{ and  }\mathcal{fD}(T_R\circ \varphi, R)=\mathcal{fD}(\varphi, R)+1.$$

{\centering
\includegraphics[width=7cm]{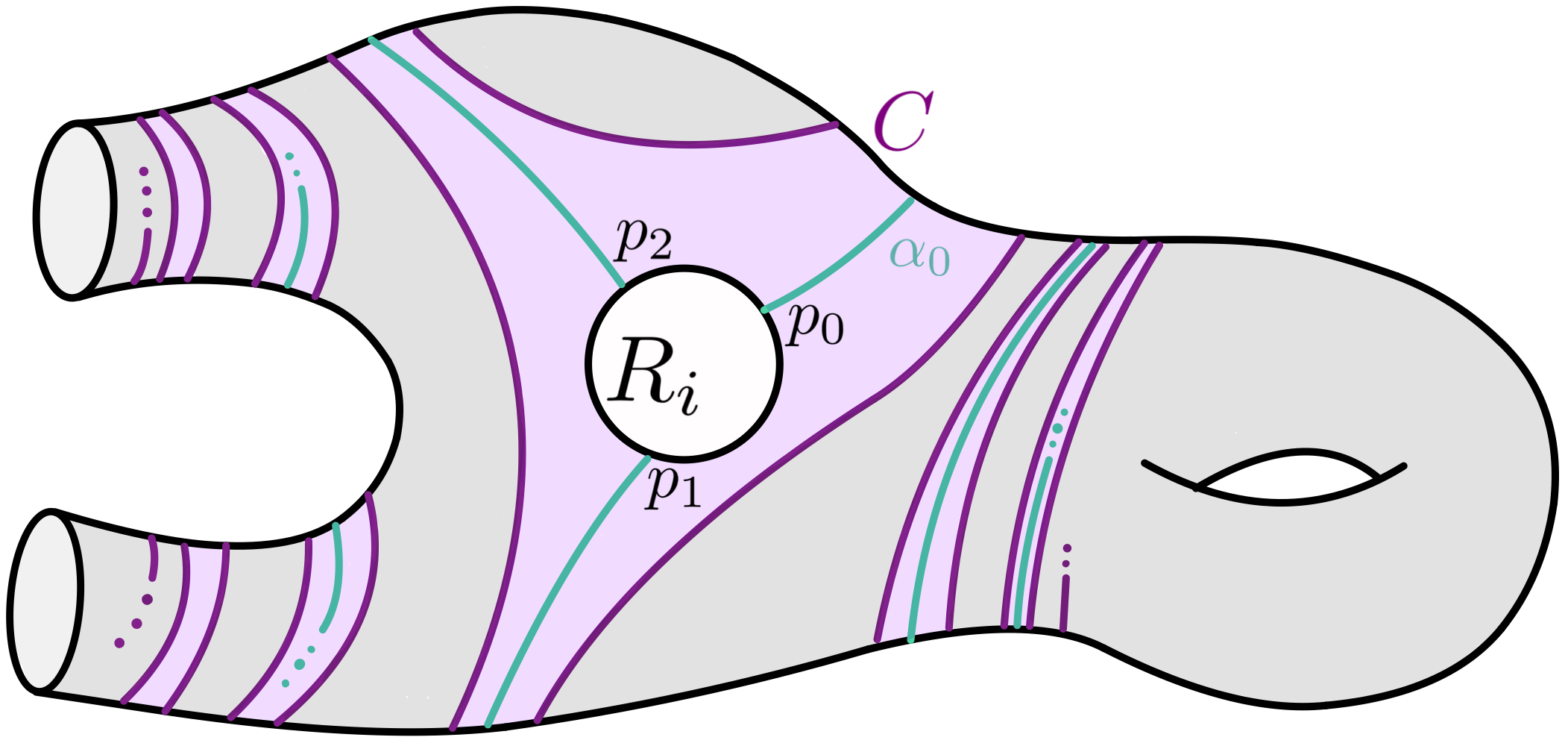}
\captionof{figure}{Surface $S$ as described in definition \ref{definition:11} with component $C \subset S \setminus L$ in purple.}
\label{figure:4}}
\vskip 0.3cm

\noindent
One can coarsely read off the fractional Dehn twist coefficient from the action of $\varphi$ on essential arcs: 

\begin{lemma}\label{lemma:read_of_FDTC} \cite[Proposition 2.9, Proposition 2.20, Proposition 2.22]{KR13} Consider a surface $S$, $\varphi\in \text{Diff}^+(S, \partial S)$ isotopic to a pseudo-Anosov map, a boundary component $R\subset \partial S$, and a simple essential arc $\gamma$ with at least one endpoint in $R$. Take a collar neighbourhood $R\times [0,1]\subset S$ and equip $S'=S\setminus R\times [0,1]$ with a hyperbolic metric with geodesic boundary. Isotope $\gamma$ relative to $R$ to the arc $\gamma_0.\gamma'.\gamma_1$ with $\gamma_0, \gamma_1\subset R\times [0,1]$ and $\hat{\gamma}\subset S'$ a geodesic meeting $\partial S'$ orthogonally. Let $\delta:=\varphi(\rho)$ and similarly define $\delta_0, \delta', \delta_1$. Then: $$|i(\gamma_0, \delta_0)-|\mathcal{fD}(\varphi, R)||\leq 1,$$
for $i(\gamma_0, \delta_0)$ denoting the unsigned (minimal) intersection number of $\gamma_0, \delta_0$ relative endpoints. 
\end{lemma}

\noindent
The fractional Dehn twist coefficient has, in particular, been studied due to its connection to the contact topology of 3-manifolds. Let us briefly review this topic.\\

\noindent
A {\it contact structure} on an oriented 3-manifold $M$ is a nowhere integrable plane field, that is a plane field locally defined by $\{\alpha=0\}$ for $\alpha$ a 2-form with $\alpha\wedge d\alpha$ a no-where vanishing volume-form agreeing with the orientation of $M$.\footnote{Some authors merely require $\alpha\wedge d\alpha\neq 0$ and call contact structures agreeing with the orientation of M positive. We will adopt the common convention that all contact structures are positive.} An embedded disc $D\subset M$ is {\it overtwisted} with respect to the contact structure $\xi$ if a neighbourhood of $D$ is homeomorphic to a neighbourhood of $D_0 := \{z = 0, \rho\leq \pi\} \subset \mathbb R^3$ via a map taking $\xi$ to $\{\cos(\rho )dz + \rho \sin(\rho)d\phi = 0\}$ (here $\rho, \phi, z$ are cylindrical coordinates on $\mathbb R^3$). Informally, an overtwisted disc is a disc $D$ with singular foliation $\xi_p\cap T_pD$ as shown in figure \ref{figure:5}. The key dichotomy in contact topology is between those contact structures admitting an overtwisted disc, called {\it overtwisted}, and all others, which are called {\it tight}. Every manifold admits an overtwisted contact structure and by \cite{Eli89} overtwisted contact structures are isotopic if and only if they are homotopic. On the other hand, tight contact structures are less well-understood, do not exist in every 3-manifold, and it is an open question whether every hyperbolic 3-manifold admits a tight contact structure.\\

{\centering
\includegraphics[width=2.8cm]{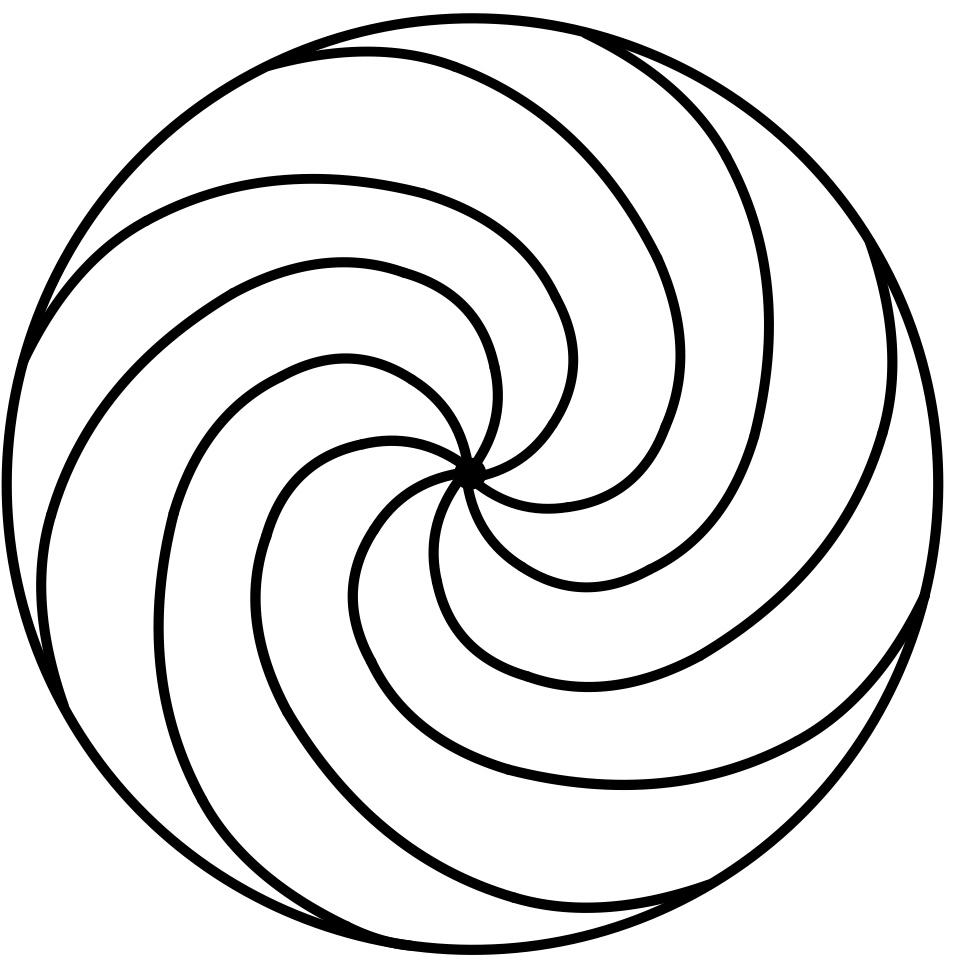}
\captionof{figure}{Sketch of the singular foliation on an overtwisted disc.}
\label{figure:5}}
\vskip 0.3cm

\noindent
An {\it open book decomposition} of a closed 3-manifold $M$ is a pair $(B, \pi)$, such that $B\subset M$ is a fibred, oriented link called the {\it binding}, and $\pi:M\setminus B\to S^1$ is a fibration of the link-complement with each fibre $\Sigma_\theta=\pi^{-1}(\theta)$ a finite type surface with boundary $\partial \Sigma_\theta=B$. We call the fibres $\Sigma_\theta$ {\it pages}. By \cite{Gir02}, see also \cite{Etn04}, there is a, so called, {\it Giroux correspondence} between contact structures on $M$, up to isotopy, and open book decompositions on $M$, up to a move called stabilisation. One may therefore ask which open book decompositions correspond to tight contact structures.\\

\noindent
Up to homeomorphism, an open book decomposition $(B, \pi)$ is determined by the corresponding {\it abstract open book decomposition} $(S, \varphi)$ where: $S$ is a fibre $\Sigma_\theta$ and $\varphi\in \text{Diff}^+(S, \partial S)$ is a monodromy of the fibration $\pi$ such that for each $R\subset \partial S$ the boundary slope $\gamma_{lat,R}$ given by $\varphi$ is the meridian of $B$.

\begin{theorem} \label{theorem:12} \cite[Theorem 4.3]{HKM08} Consider an open book decomposition $(B, \pi)$ for the closed manifold $M$ with corresponding abstract open book decomposition $(S, \varphi)$. If $M\setminus B$ is hyperbolic, the binding is connected, and $\mathcal{fD}(\varphi)\geq 1$, then the contact structure corresponding to $(B, \pi)$ is tight. 
\end{theorem}

\section{Pointed Arc Graph}\label{section:5} 
\noindent 
To understand the cusp-skew $\mathcal{sk}(\varphi, R)$, we introduce a modified version of the arc graph.

\begin{definition}\label{definition:13}
Consider a compact surface $S$ with boundary and $\chi(S)<0$, choose a component $R\subset \partial S$, and a point $p\in R$. We take {\it essential simple arc} to mean a smooth proper map $[0, 1]\to S$, that is an embedding on its interior and is not homotopic through proper maps to a point. We define the {\it pointed arc graph} $\mathcal{A}_p(S)$ as follows: The vertices of $\mathcal{A}_p=\mathcal{A}_p(S)$ are essential simple arcs $\gamma$ in $S$ with $\gamma\cap R=\{p\}$, modulo isotopy of $S$ fixing $p$. Vertices are adjacent if they have representatives that are disjoint away from $p$.
\end{definition}

\noindent
Note, that in the above definition an essential arc must not be homotopic into a puncture even with respect to a homotopy not fixing $\partial S$. Moreover, this definition requires that arcs of $\mathcal{A}_p(S)$ have at least one endpoint in $R$ and meet $R$ only at $p$. Given an arc $\gamma$ in $S$ with at least one endpoint at $p$, let $[\gamma]_p\in \mathcal{A}_p(S)$ and $[\gamma]\in \mathcal{A}(S,R)$ denote its based and unbased isotopy classes, respectively. Let $\mathcal{A}^1(S,R), \mathcal{A}^2(S,R)\subset \mathcal{A}(S,R)$ denote the sets of arcs with one or two endpoints in $R$ respectively. 
We now describe an identification of the vertex set of $\mathcal{A}_p(S)$ and the set $(\mathcal{A}^1(S,R)\times \mathbb Z)\cup (\mathcal{A}^2(S,R)\times \mathbb Z[\frac 12])\text{, for }$ $\mathbb Z[\frac 12]:=\{0, \frac 12\}+\mathbb Z$. We will then give a coarse description of the edge set of $\mathcal{A}_p(S)$ in terms of this identification.\\

\noindent
Consider an arc $\gamma$ representing a vertex $[\gamma] \in\mathcal{A}(S,R)$ with $\gamma\cap R=\{p\}$. Take a collar neighbourhood $R \times [0, 1] \subset S$, let $S' = S \setminus (R \times I)$, and endow $S'$ with a hyperbolic metric with geodesic boundary. Then $\gamma$ is isotopic relative $p$ to an arc $\hat\gamma$ which is geodesic in $S'$, orthogonal to $\partial S'$, and intersects $R \times I$ in arcs disjoint away from $p$. The geodesic part $\hat\gamma \cap S'$ is determined by the free isotopy class $[\gamma] \in\mathcal{A}(S, R).$\\

{\centering
\includegraphics[width=7.5cm]{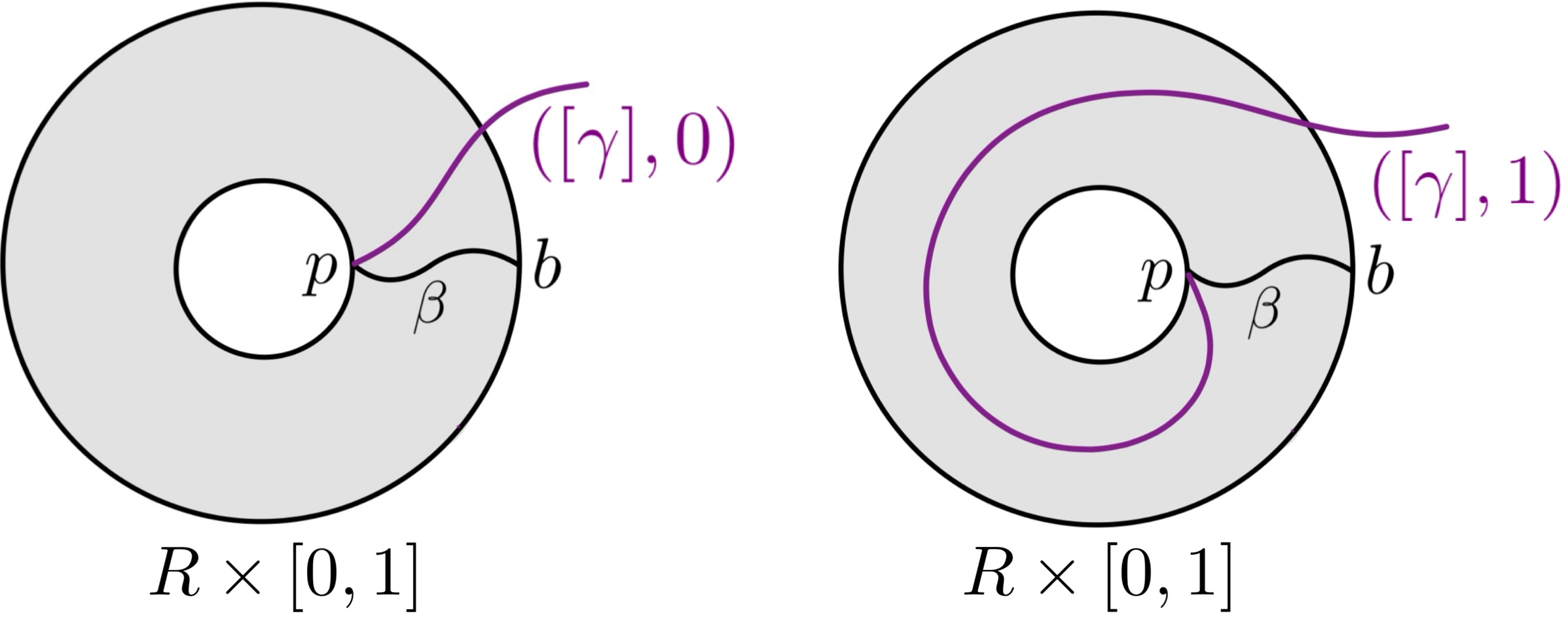}
\captionof{figure}{$R \times [0, 1]$ with arcs $([\gamma], 0)$ and $([\gamma],1)$ in purple.}
\label{figure:6}}
\vskip 0.3cm

\noindent
There are countably many isotopy classes of arcs in $S'$, so we may choose $b \in R \times \{1\} \subset  \partial S'$ which is not the endpoint of any geodesic arc $\hat\gamma\cap S'$ and we may choose an arc $\beta \in R\times[0,1]$ from $p$ to $b$. Suppose $\gamma$ has one endpoint in $R$. As shown in figure \ref{figure:6}, there are two isotopy classes of simple arcs in $R \times [0, 1] - \beta$ connecting $p$ to the endpoint of the geodesic $\hat\gamma\cap S'$. Hence there are two choices, up to isotopy relative $p$, for $\hat\gamma$ in $S - \beta$, which we denote $([\gamma], 0)$ and $([\gamma], 1)$ as shown in figure \ref{figure:6}. Notice that $([\gamma], 1) = T_R([\gamma], 0)$, where $T_R$ denotes a Dehn twist about $R$. All possibilities for $\hat\gamma$, up to isotopy fixing $R$, can be obtained this way for suitable choice of $\beta$ and all choices for $\beta$ are related by Dehn twist $T_R$. Hence all $[\gamma]_p \in \mathcal{A}_p(S)$ with fixed free isotopy class $[\gamma] \in \mathcal{A}^1(S, R)$ are of the form $T_R^n([\gamma], 0), n \in \mathbb Z$.\\

{\centering
\includegraphics[width=11.3cm]{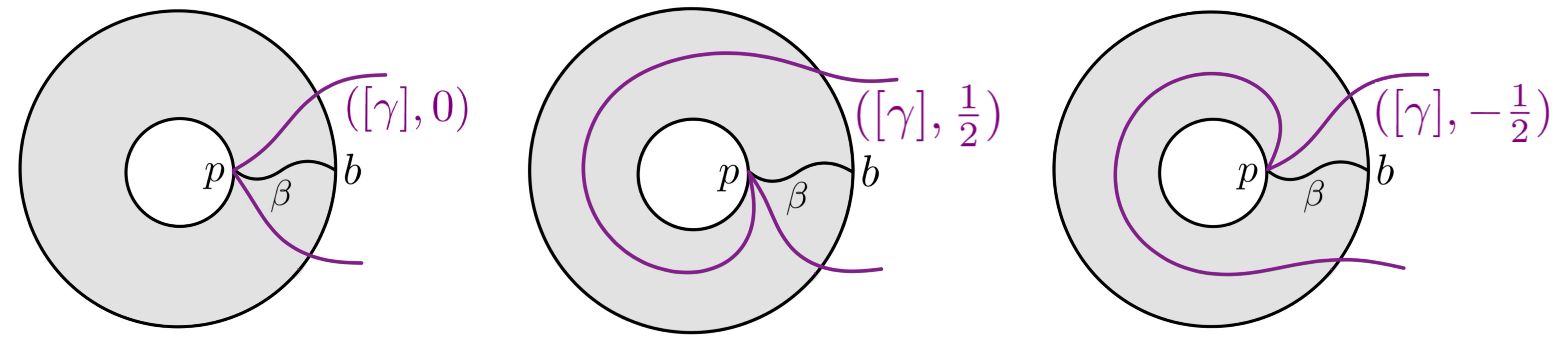}
\captionof{figure}{$R \times [0, 1]$ with arcs $([\gamma], 0), ([\gamma], \frac 12)$, and $([\gamma],-\frac 12)$ in purple.}
\label{figure:7}}
\vskip 0.3cm

\noindent
Now suppose $\gamma$ has two endpoints in $R$. As shown in figure \ref{figure:7}, there are three isotopy classes of pairs of disjoint simple arcs in $R \times [0, 1] - \beta$ connecting $p$ to the endpoints of $\hat\gamma \cap S'$. This gives three choices, up to isotopy relative $p$, for $\hat \gamma$ in $S-\beta$, which we denote $([\gamma],-\frac 12), ([\gamma],0),$ and $([\gamma],\frac 12)$ as shown in figure \ref{figure:7}. Observe that, $([\gamma],-1) = T_R^{-1}([\gamma], 1)$. Hence all $[\gamma]_p \in \mathcal{A}_p(S)$ with fixed free isotopy class $[\gamma] \in \mathcal{A}^2(S,R)$ are of the form $T_R^n([\gamma], 0),$ or $T_R^n([\gamma], \frac 12 ), n \in \mathbb Z.$ We introduce the notation: $([\gamma], n) := T_R^n([\gamma], 0)$, and $([\gamma], n + \frac 12 ) := T_R^n([\gamma], \frac 12 )$. In summary, we have described a surjection from $(\mathcal{A}^1(S,R)\times\mathbb Z)\cup(\mathcal{A}^2(S,R)\times\mathbb Z[\frac 12])$ to the vertex set of $\mathcal{A}_p(S)$. This map depends on our choice of $b, \beta,$ and metric on $S'$. Unless stated otherwise, all results will hold for any such choices. We now coarsely describe how adjacency in $\mathcal{A}_p(S)$ behaves under this surjection.

\begin{lemma}\label{lemma:14} Consider the above setup and take $[\delta],[\gamma] \in \mathcal{A}(S,R), n,m \in \mathbb Z[\frac 12]$. Under the above surjection:
\begin{enumerate}
\item for $[\gamma] \in \mathcal{A}^1(S, R), ([\gamma], n), ([\gamma], m) \in \mathcal{A}_p(S)$ are adjacent if and only if $|m - n| \leq 1$;
\item for $[\gamma] \in \mathcal{A}^2(S, R), ([\gamma], n), ([\gamma], m) \in \mathcal{A}_p(S)$ are adjacent if and only if $|m - n| < 1$;
\item if $([\gamma], n), ([\delta], m) \in \mathcal{A}_p(S)$ are adjacent, then $[\gamma], [\delta] \in \mathcal{A}(S, R)$ are adjacent and $|m - n| \leq 1$;
\item consider the above discussed representative $\hat \gamma$ of $([\gamma], n) \in \mathcal{A}_p(S)$ and denote the initial segment of $\hat\gamma \cap (R \times I)$ by $\gamma_0$. Similarly define $\delta_0$ for $([\delta], m) \in \mathcal{A}_p(S)$, with $[\delta]\neq [\gamma]$. Let $i(\gamma_0, \delta_0)$ denote the unsigned intersection number of $\gamma_0, \delta_0$ relative endpoints. Then $||m - n| - i(\gamma_0, \delta_0)| \leq 2$.
\end{enumerate}
\end{lemma}

\noindent
{\it Proof:} We define a {\it bigon} between two arcs $\gamma, \delta$ in $S$ to be an embedded disc in $S \setminus (\gamma \cup \delta)$ whose boundary is a union of two closed arcs $a \subsetneq \gamma, b \subsetneq \delta$. It is a standard fact, that a pair of transverse arcs $\gamma, \delta$ realise minimal intersection numbers in their isotopy classes if and only if they have no bigons. Hence, for each statement it is sufficient to find suitable bigon-free representatives of the relevant isotopy classes. Recall that when we describe an element of $\mathcal{A}_p(S)$ as $([\gamma],n)$ for $[\gamma] \in \mathcal{A}(S,R)$, we fixed the following conventions: $R \times [0, 1]$ denotes a collar neighbourhood of $R, S' = S - R \times [0, 1]$, and $\hat \gamma$ is a representative of $([\gamma],n)$ which is geodesic on $S'$. Similarly define $\hat\delta$.\\

\noindent
(1): By swapping arcs, we may assume that $m \geq n$. By Dehn twisting in $R$, we may assume that $n = 0$. In figure \ref{figure:8} left, we see representatives of $([\gamma],0)$ and $([\gamma],1)$ disjoint away from $p$. In figure \ref{figure:8} right, we see intersecting, bigon-free representatives of $([\gamma],0)$ and $([\gamma],2)$. It is clear how to generalise figure \ref{figure:8} to find intersecting, bigon-free representatives of $([\gamma], 0)$ and $([\gamma], m)$ for $m > 2$. This concludes (1).\\

{\centering
\includegraphics[width=7.3cm]{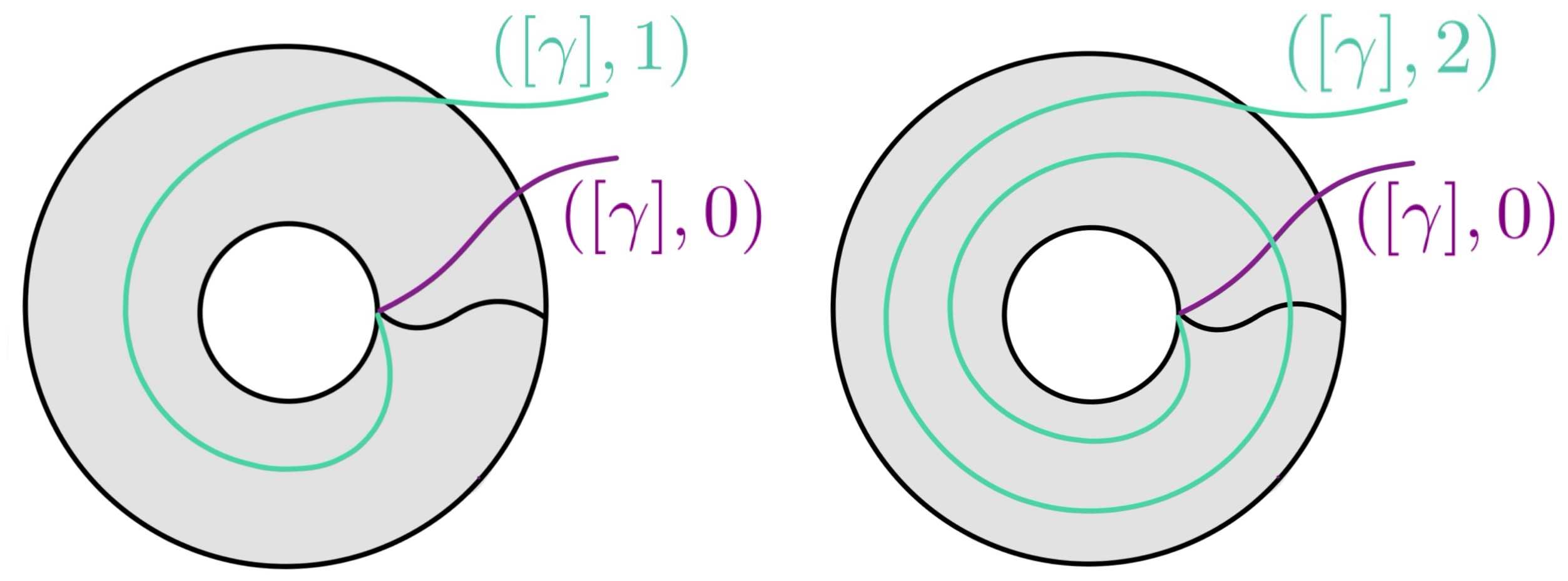}
\captionof{figure}{Bigon-free representative of $([\gamma], 0), ([\gamma], 1)$, and $([\gamma], 2)$.}
\label{figure:8}}
\vskip 0.3cm

\noindent
(2): As in (1), we may assume $m > n$ and $n = 0$ or $\frac 12$. First suppose $n = 0$. In figure \ref{figure:9} left we see bigon-free representatives of $([\gamma], 0)$ and $([\gamma], \frac 12 )$ and of $([\gamma], 0)$ and $([\gamma], 1)$. It is clear how to generalise figure \ref{figure:9} right to show $([\gamma], 0)$ and $([\gamma], m)$ are adjacent if and only if $m \leq 1$.\\

{\centering
\includegraphics[width=8cm]{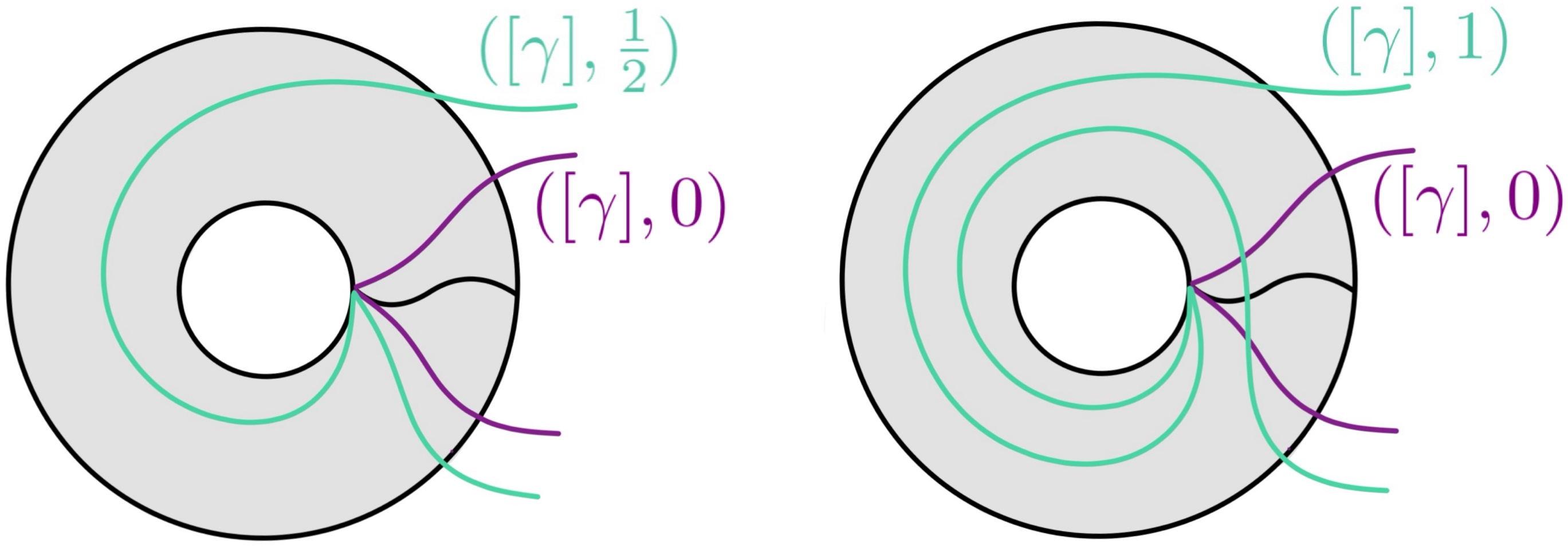}
\captionof{figure}{Bigon-free representative of $([\gamma], 0), ([\gamma], \frac 12 ),$ and $([\gamma], 1)$.}
\label{figure:9}}
\vskip 0.3cm

\noindent
Now suppose $n = \frac 12$ . In figure \ref{figure:10} left we see representatives of $([\gamma], \frac 12 )$ and $([\gamma], 1)$ disjoint away from $p$. On the other hand, as we see in figure \ref{figure:10} right and its generalisations, $([\gamma], \frac 12 )$ and $([\gamma], m)$ are non-adjacent for $m > 1$. This concludes (2).\\

{\centering
\includegraphics[width=8cm]{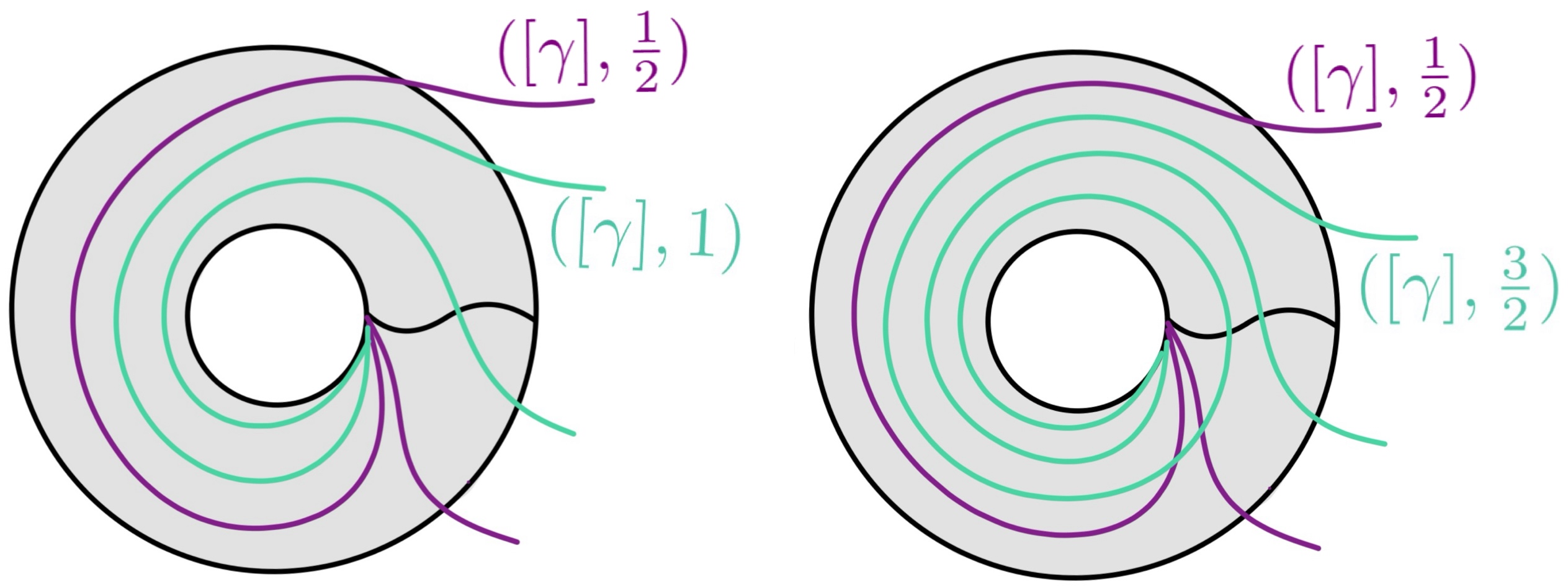}
\captionof{figure}{Bigon-free representative of $([\gamma], \frac 12 ), ([\gamma], 1),$ and $([\gamma], \frac 32)$.}
\label{figure:10}}
\vskip 0.3cm

\noindent
(3): Certainly for $([\gamma], n)$ and $([\delta], m) \in \mathcal{A}_p(S)$ to have representatives disjoint away from $p$, $[\gamma], [\delta] \in \mathcal{A}(S, R)$ must have disjoint representatives, so assume this is the case. By possibly twisting in $R$, we may assume that $n = 0$ or $n = -\frac 12$. The elements $([\gamma],n),n = 0$ or $-\frac 12$ and $([\delta],m), m = \frac 32$ or $m = 2$ have representatives as shown in figure \ref{figure:11}. Observe that the case shown in figure \ref{figure:11} left is bigon-free. If the blue-shaded area in figure \ref{figure:11} right is a bigon, we may isotope it away and be in the setup of figure \ref{figure:11} left. Similarly, we see that $([\gamma], n), n = 0$ or $-\frac 12$ are not adjacent to $([\delta],m),m>2$.\\

{\centering
\includegraphics[width=10cm]{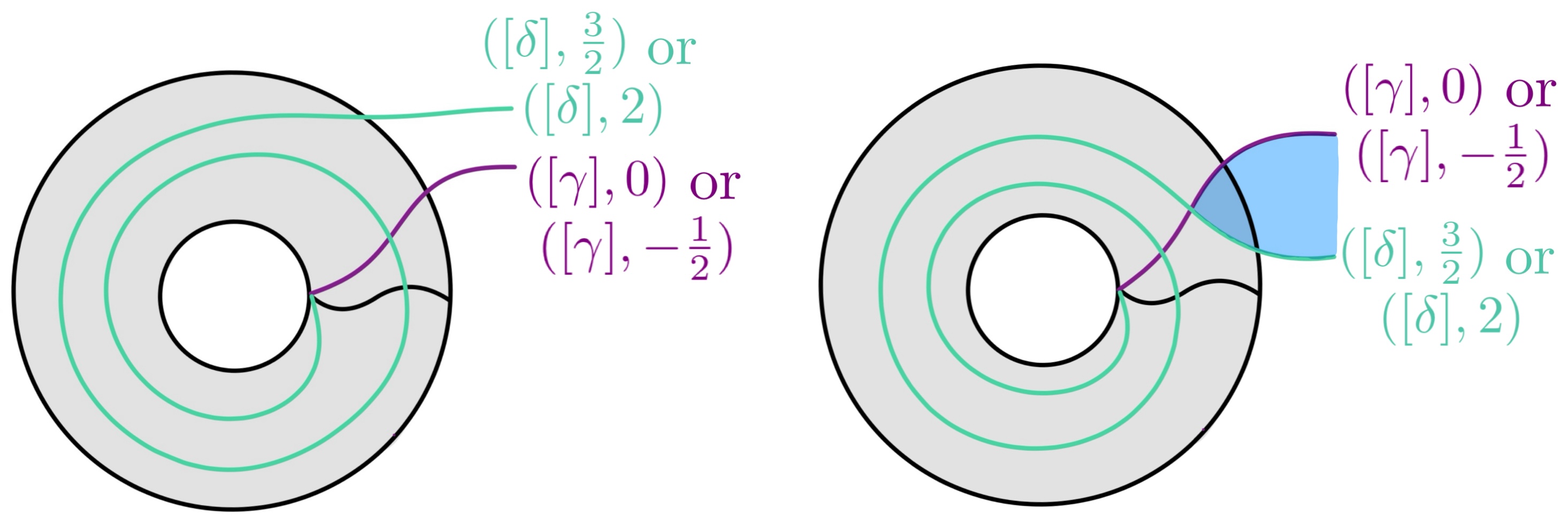}
\captionof{figure}{Representatives of the initial segments of $([\gamma], n)$ and $([\delta], m)$ for $n = 0$ or $-\frac 12$ and $m = \frac32$ or $2$.}
\label{figure:11}}
\vskip 0.3cm

\noindent
It remains to show that $([\gamma], - \frac 12 )$ and $([\delta], 1)$ are not adjacent. These have representatives as shown in figure \ref{figure:12}. Again, note that figure \ref{figure:12} left is bigon-free and figure 12 right is either bigon-free or can be isotoped to the case of figure \ref{figure:12} left. This concludes (3).\\

{\centering
\includegraphics[width=8cm]{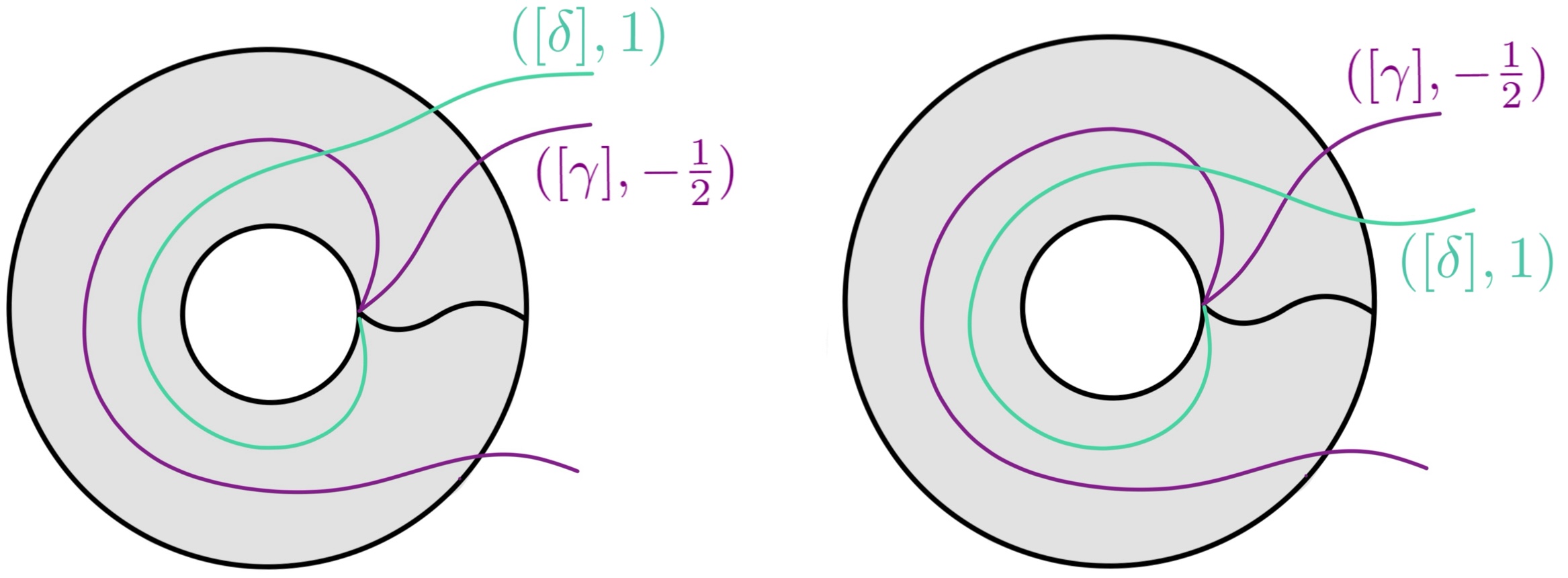}
\captionof{figure}{Representatives of $([\gamma], -\frac 12 )$ and $([\delta], 1)$.}.
\label{figure:12}}
\vskip 0.3cm

\noindent
(4): Recall that $\hat \gamma \cap S', \hat \delta \cap S'$ are geodesics and therefore have distinct endpoints. By Dehn twisting both arcs in R, we may assume $\gamma_0$ is as shown in figure \ref{figure:13}. Hence $n = 0,-\frac 12$, or $-1$. There are two cases for the first endpoint of $\hat \delta \cap S'$ as shown in figure \ref{figure:13}. Suppose $\delta_0$ wraps $k$-times anti-clockwise around $R$. Then $m$ is one of $k,k- \frac 12,$ or $k-1$. Moreover $i(\gamma_0,\delta_0)$ is $k-1$ or $k$ depending on whether we are in case 1 or 2. In summary, $k-1\leq m-n\leq k+1$ and $||m-n|-i(\gamma_0,\delta_0)|\leq 2$. This concludes (4). $\square$ \\

{\centering
\includegraphics[width=8cm]{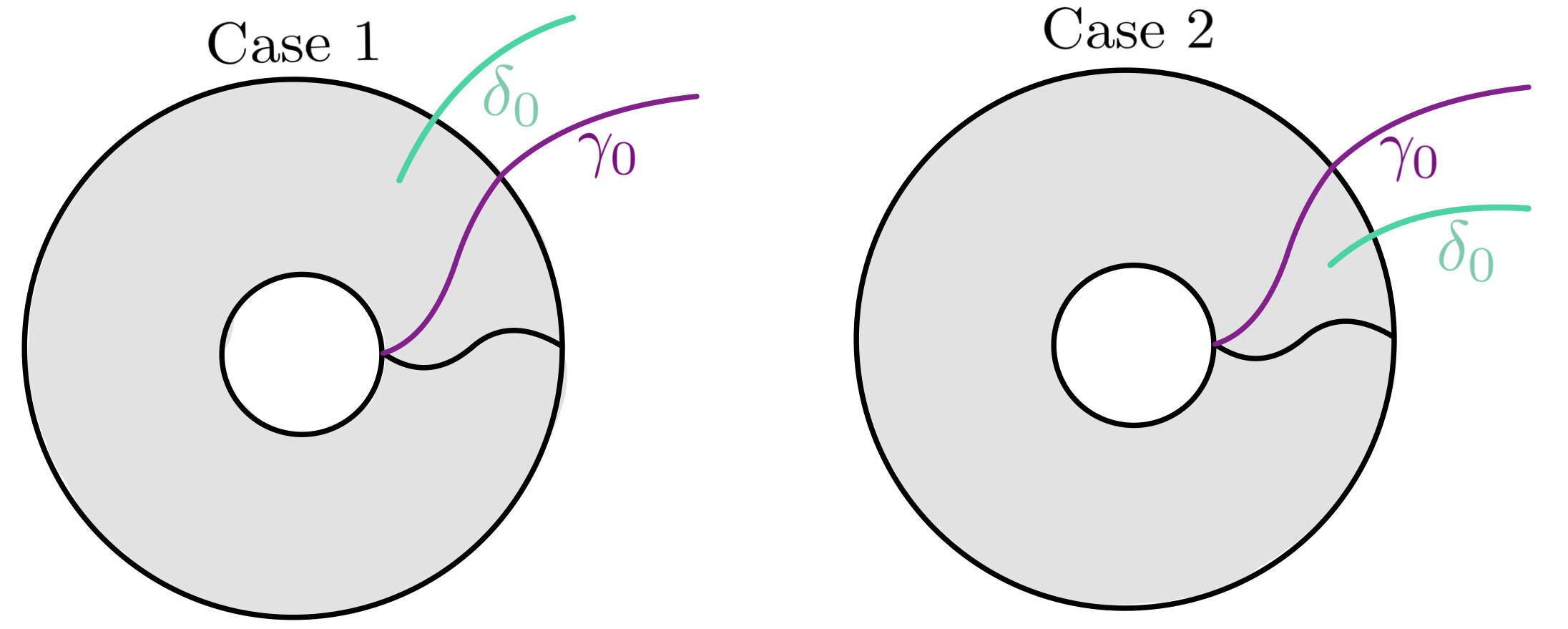}
\captionof{figure}{Two cases of relative position of $\gamma_0$ and $\delta_0$ in the proof of lemma \ref{lemma:14} (4).}
\label{figure:13}}
\vskip 0.3cm

\begin{corollary} \label{corollary:15} The above map from $(\mathcal{A}^1(S,R)\times \mathbb Z)\cup (\mathcal{A}^2(S,R)\times \mathbb Z[\frac 12])$ to the vertices of $\mathcal{A}_p(S)$ is a bijection. 
\end{corollary}

\noindent
{\it Proof:} From the above discussion, we know this map is a surjection. Suppose $([\gamma],n)$ and $([\delta],m)$ are equal in $\mathcal{A}_p(S)$. Then certainly $[\gamma], [\delta]$ must be equal in $\mathcal{A}(S, R)$. By lemma \ref{lemma:14}, $([\gamma], n)$ and $([\gamma], m)$ have the same neighbours in $\mathcal{A}_p(S)$ if and only if $n = m$. In particular, for $n\neq m$: $([\gamma], n)$ and $([\gamma], m)$ are not isotopic relative $p$ and the discussed map is injective. $\square$

\begin{corollary} \label{corollary:16} The pointed arc graph $\mathcal{A}_p(S)$ is 2-quasi-isometric to the graph $G$ with vertex set $(\mathcal{A}^1(S, R) \times \mathbb Z)\cup(\mathcal{A}^2(S,R)\times\mathbb Z[\frac 12])$ and an edge between $([\gamma],n)$ and $([\delta],m)$ if $n = m$ and $[\gamma],[\delta]$ are adjacent in $\mathcal{A}(S,R)$ or if $[\gamma] = [\delta]$ and $|m - n| \leq 1$.
\end{corollary}

\noindent
{\it Proof:} If $[\gamma]\in \mathcal{A}(S,R)$ and $|m-n|\leq 1$, then by lemma \ref{lemma:14}: $([\gamma],n)$ and $([\gamma],m)$ are adjacent in $\mathcal{A}_p(S)$ or have a mutual neighbour $([\gamma], \frac{m+n}2)\in \mathcal{A}_p(S)$. If $[\gamma], [\delta]\in \mathcal{A}(S,R)$ are adjacent, then $([\gamma], n), ([\delta],n)$ are adjacent in $\mathcal{A}_p(S)$. Now suppose that $([\gamma], n)$ and $([\delta],m)\in \mathcal{A}_p(S)$ are adjacent. Then by lemma \ref{lemma:14}: $[\gamma], [\delta]\in \mathcal{A}(S,R)$ are equal or adjacent and $|m-n|\leq 1$. Hence $([\gamma], n), ([\delta],n), ([\delta],m)$ is a path of length at most $2$ in $G$. $\square$\\

\noindent
Finally, this identification of the vertex set of $\mathcal{A}_p(S)$ allows us to coarsely read off the fractional Dehn twist coefficient of $\varphi \in \text{Diff}^+(S, R)$.

\begin{lemma} \label{lemma:17}
Consider a surface $S$ with $\chi(S) < 0$, distinguished boundary component $R \subset \partial S, p \in R$, and a $\varphi \in \text{Diff}^+(S,\partial S)$ isotopic to a pseudo-Anosov map. Suppose that under an identification of $\mathcal{A}_p(S)$ with $(\mathcal{A}^1(S,R)\times\mathbb Z)\cup(\mathcal{A}^2(S,R)\times\mathbb Z[\frac 12])$, as described above, $\varphi$ maps $([\delta],n)$ to $([\rho],m)$. Then: 
$$|\mathcal{fD}(\varphi, R)-(m-n)|\leq 5.$$
\end{lemma}

\noindent
{\it Proof:} Recall that, as we see in figure \ref{figure:14}, when identifying $([\delta],n)$ with an isotopy class in $\mathcal{A}_p(S)$, we take a representative $\hat{\delta} = \delta_0.\delta'.\delta_1$, with $\delta_0, \delta_1$ arcs in the collar neighbourhood $R \times [0, 1]$ and $\delta'$ a geodesic in $S' =S-(R\times[0,1])$. Similarly, we obtain a representative $\hat{\rho}=\rho_0.\rho'.\rho_1$ of $([\rho], m)$. Since composing $\varphi$ with a Dehn twist $T_R$ in $R$ increases both $m - n$ and $\mathcal{fD}(\varphi, R)$ by 1, we may assume that $|m - n| \leq 1$. Thus by lemma \ref{lemma:14}, $i(\delta_0, \rho_0) \leq 3$. By lemma \ref{lemma:read_of_FDTC}: $||\mathcal{fD}(\varphi,R)|-i(\delta_0,\rho_0)|\leq1$. In summary, $|\mathcal{fD}(\varphi,R)-(m-n)|\leq5$. $\square$\\

{\centering
\includegraphics[width=5.5cm]{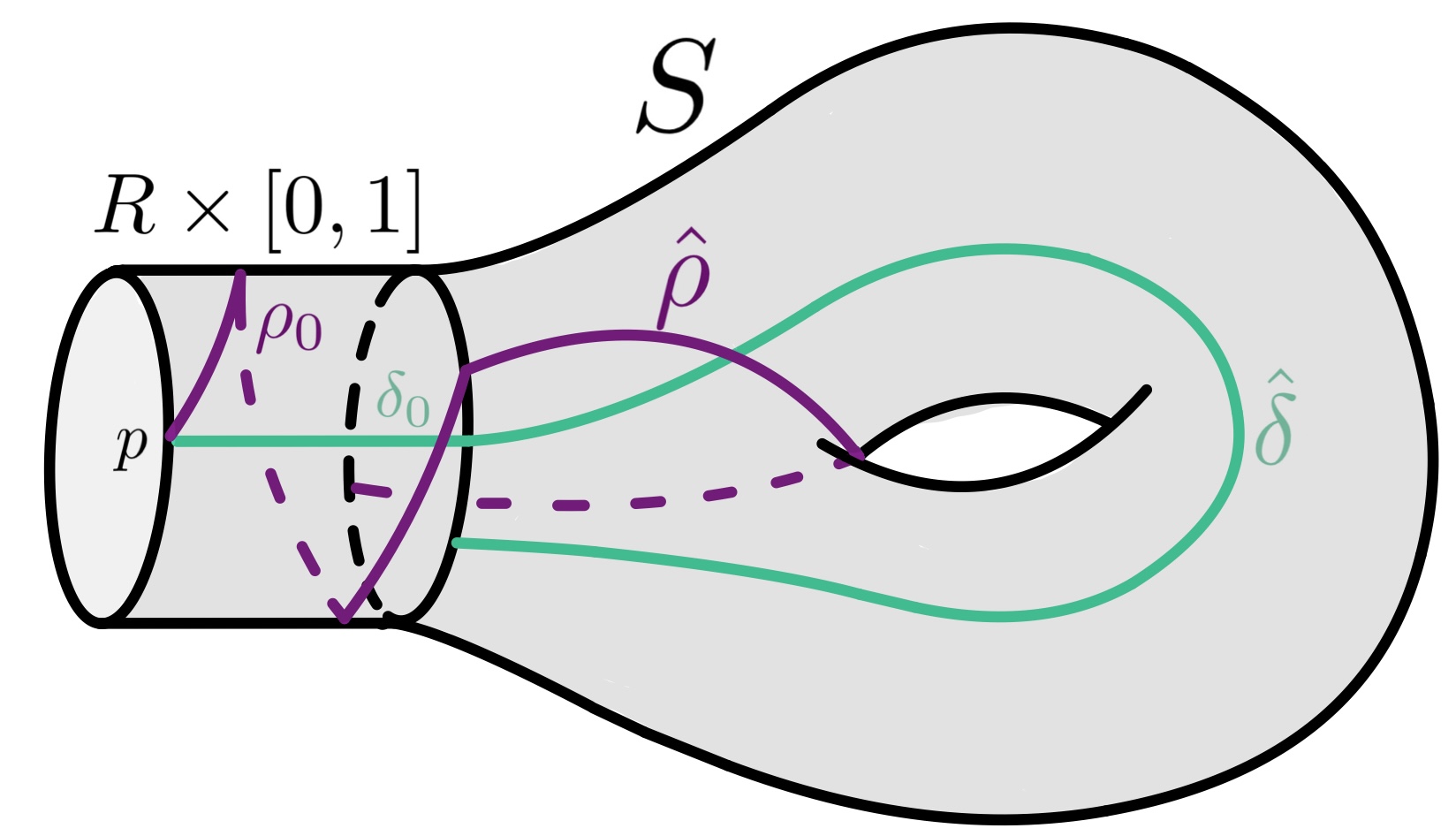}
\captionof{figure}{Sketch of $\hat{\delta}$ and $\hat\rho$ on $S$.}
\label{figure:14}}
\vskip 0.3cm

\section{Proof of Theorem \ref{theorem:2}} \label{section:6}
\noindent
Consider again the setup from the introduction, that is $S$ is an oriented surface with distinguished boundary $R \subset \partial S, \varphi  \in \text{Diff}^+ (S, \partial S)$ is isotopic to a pseudo-Anosov diffeomorphisms, and $T = R \times  S^1$ is a boundary torus of the mapping torus $M_\varphi $. We may define the translation distance $d_{\mathcal{A}_p} (\varphi) := \inf_{\rho \in \mathcal{A}_p(S)} d(\rho , \varphi (\rho ))$.

\begin{lemma} \label{lemma:18} Consider the above setup and recall $\gamma _{lat,R} \in  \pi _1(T)$ is given by $\{q\} \times  S^1$. Then:
$$l(\gamma_{lat,R})\leq 6|\chi(S)|\cdot d_{\mathcal{A}_p}(\varphi).$$
\end{lemma}

\noindent
{\it Proof:} Let $N = d_{\mathcal{A}_p}(\varphi )$ and take a path of arcs $\rho_0,...,\rho_N = \varphi(\rho_0): [0,2] \to S$ in $\mathcal{A}_p(S)$ such that $\rho_i(0) = p$. Up to isotopy relative $\partial S$, we may assume that $\varphi $ pointwise fixes a collar neighbourhood $R \times  [0, 1] \subset S$ of $R$. Moreover, we may assume that $\rho_i(t) = (p, t) \in  R \times  [0, 1]$ for $t \in  [0, 1]$. Consider the universal cover $\unicov{M_\varphi}\to M_\varphi$ and let $\unicov{p}$ be a lift of $\,\quotient{((p,\frac 12)\times {0})}{\sim} \in M_\varphi$. Consider lifts $\unicov{\rho_i}\subset \unicov{M_\varphi}$ at $\unicov{p}$ of $\,\quotient{(\rho_i\times {0})}{\sim} \subset M_\varphi $ and let $\unicov{g_i}$ denote their respective geodesic representatives in $\unicov{M_\varphi}$. Finally, let $C$ be the maximal cusp of $M_\varphi$ corresponding to $T$, let $\unicov{C}\subset M_\varphi$ be the horoball lift of $C$ at $\unicov{p}$, and let $\unicov{q_i}:=\unicov{g_i}\cap\partial \unicov{C}$.\\

\noindent
Notice that $H(s,t) := \,\quotient{(\rho_0(s)\times {t})}{\sim}, t \in  [0,1]$ is an isotopy between $\,\quotient{(\rho_0 \times {0})}{\sim} $ and $\quotient{(\rho_N \times {0})}{\sim} $ and $H(\frac 12,t)$ is homotopic to $\gamma_{lat,R}$ . Thus $\unicov{\rho_0}, \unicov{\rho_N}$ differ by the deck transform given by $\gamma_{lat,R} \in \pi_1(M_\varphi)$ and therefore a path from $\unicov{q_0}$ $\unicov{q_N}$ in $\unicov{C}$ descends to a loop representing $\gamma_{lat,R} \in  \pi_1(T )$. It is thus sufficient to show that in the induced Euclidean metric on $\partial C$, we have $d(\unicov{q_i}, \unicov{q_{i+1}}) \leq  6|\chi(S)|$.\\

\noindent
Suppose first that the unbased isotopy classes $[\rho_i],[\rho_{i+1}] \in  \mathcal{A}(S,R)$ are equal. Then under an identification of $\mathcal{A}_p(S)$ with $(\mathcal{A}^1(S,R)\times \mathbb Z )\cup(\mathcal{A}^2(S,R)\times \mathbb Z [\frac 12])$, as discussed in section \ref{section:5}, $[\rho_i],[\rho_{i+1}]$ correspond to $([\rho ],n), ([\rho ],m)$ with $|m-n|\leq 1$. In particular, $\unicov{q_i}$ and $\unicov{q_{i+1}}$ are equal or differ by a lift of $\gamma_{\partial,R}$. So in $\partial \unicov{C}$, we have $d(\unicov{q_i}, \unicov{q_{i+1}}) \leq  6|\chi(S)|$.\\

\noindent
Thus we may assume that $\rho_i, \rho_{i+1}$ extend to an ideal triangulation $\tau$ of $S$. Consider the pleating $f : S \to M_\varphi $ of the embedding $S \to S \times  \{0\}/\sim $ along $\tau$. As in lemma \ref{lemma:9}, take a sufficiently small cusp $C_0 \subset  C$ such that $f(S)$ intersects $C_0$ in tips of ideal triangles. Denote the distance between $C_0$ and $C$ by $d$, the lift of $C_0$ within $\unicov{C}$ by $\unicov{C_0}$, and the curve $f(S) \cap \partial C_0$ by $\gamma_0$. Just as in the proof of lemma \ref{lemma:9}, we know that $\text{length}(\gamma _0) \leq  e^{-d}\cdot  6|\chi(S)|.$ By considering a lift $\unicov{f}:\unicov{S} \to \unicov{M_\varphi}$ of $f$ containing $\unicov{g_i}$ and $\unicov{g_{i+1}}$, we see that $\unicov{g_i} \cap \unicov{C_0}$ and $\unicov{g_{i+1}} \cap \unicov{C_0}$ are contained in a lift $\unicov{\gamma_0}$ of $\gamma_0$. Hence in the induced Euclidean metric on $\partial \unicov{C_0}: d(\unicov{g_i} \cap \unicov{C_0}, \unicov{g_{i+1}} \cap \unicov{C_0}) \leq  e^{-d}\cdot  6|\chi(S)|$. Since the Euclidean metrics on $C_0$ and $C$ are related by a dilation of factor $e^d$, we have the desired bound on $d(\unicov{q_i}, \unicov{q_{i+1}})$ and complete the proof. $\square$\\

\noindent
Let us now recall and prove theorem \ref{theorem:2}.\\

\noindent
{\bf Theorem \ref{theorem:2}.} Consider an oriented, compact surface $S, \varphi  \in  \text{Diff}^+ (S, \partial S)$ isotopic to a pseudo-Anosov map, a boundary component $R \subset  \partial S$, and the corresponding boundary component $T = R \times  S^1$ of the mapping torus $M_\varphi $. Then we have the following inequality:
$$|\mathcal{fD}(\varphi, R)-\mathcal{sk}(\varphi, R)|\leq 6|\chi(S)|\cdot d_{\mathcal{A}(S,R)}(\varphi)+3.$$

\noindent
{\it Proof:} Observe that both sides of this inequality are unchanged when composing $\varphi $ with a Dehn twist $T_R$. Therefore, without loss of generality assume that $|\mathcal{fD}(\varphi ,R)| \leq  1$. Take $k \in  N$ and suppose $\varphi^k$ maps $(\delta,0)\in \mathcal{A}_p(S)$ to $(\rho ,n)$. By lemma \ref{lemma:17}, we see that $|n|\leq |\mathcal{fD}(\varphi^k)|+5\leq k+5$. By lemma \ref{lemma:18}, we have that:
$$l(\gamma _{lat, R}(\varphi )) = \frac{l(\gamma _{lat, R}(\varphi ^k))}{k} \leq  \frac{6|\chi(S)| }{k}\cdot d_{\mathcal{A}_p} ((\delta, 0), (\rho , n)).$$
\noindent
Suppose that $n \geq 0$ and $\rho $ has one endpoint in $R$, that is $\rho  \in  \mathcal{A}^1(S, R)$. Then by lemma \ref{lemma:14}:
$$d_{\mathcal{A}_p} ((\rho , 0), (\rho , n)) \leq  d_{\mathcal{A}_p} ((\rho , 0), (\rho , 1)) + . . . + d_{\mathcal{A}_p} ((\rho , n - 1), (\rho , n)) = |n|.$$
\noindent
If instead $\rho  \in  \mathcal{A}^2(S, R)$, then:
$$d_{\mathcal{A}_p} ((\rho , 0), (\rho , n)) \leq  d_{\mathcal{A}_p} ((\rho , 0), (\rho , 1/2)) + . . . + d_{\mathcal{A}_p} ((\rho , n - 1/2), (\rho , n)) = 2|n|.$$
Therefore, we obtain the following inequality:
$$l(\gamma _{lat, R}(\varphi))\leq  \frac{6|\chi(S)|}{k}\cdot \left(d_{\mathcal{A}_p}((\delta, 0), (\rho,0))+2|n|\right)\leq 6|\chi(S)|\cdot \frac{d_{\mathcal{A}(S,R)}(\delta,\varphi^k(\delta))+2k+10}{k}.$$
Taking the limit $k \to \infty$ gives $l(\gamma_{lat,R}(\varphi )) \leq  6|\chi(S)|\cdot  d_{\mathcal{A}(S,R)}(\varphi ) + 2$. From the definition, it is clear that $|\mathcal{sk}(\varphi ,R)| \leq  l(\gamma_{lat,R}(\varphi ))$. Finally, recall that $|\mathcal{fD}(\varphi ,R)| \leq  1$. This completes the proof. $\square$

\section{Applications to Hyperbolic Geometry}\label{section:7} 

\noindent
Theorem \ref{theorem:2} relates the cusp geometry of a fibred hyperbolic manifold to the, in some setups, better understood quantities $\mathcal{fD}(\varphi, R)$ and $\mathcal{A}(S, R)$. Let us collect some applications of this result.\\

\noindent
The fact that $\mathcal{fD}$ is a quasi-morphism, that is $|\mathcal{fD}(\varphi \circ \psi, R) - \mathcal{fD}(\varphi, R) - \mathcal{fD}(\psi, R)| \leq 1$, gives the following additivity result for cusp shape:

\begin{corollary} \label{corollary:19} Given a compact surface $S$ with distinguished boundary component $R \subset \partial S$ and maps $\varphi, \psi \in  \text{Diff}^+(S, \partial S)$ isotopic to pseudo-Anosov maps, we have:
$$|\mathcal{sk}(\varphi \circ \psi, R) - \mathcal{sk}(\varphi, R) - \mathcal{sk}(\psi, R)| \leq 6|\chi(S)| \cdot (d_{\mathcal{A}(S,R)}(\varphi) + d_{\mathcal{A}(S,R)}(\psi) + d_{\mathcal{A}(S,R)}(\varphi \circ \psi)) + 10.$$
\end{corollary}

\noindent
Given a surface $S$ and $\varphi \in  \text{Diff}^+(S, \partial S)$ we may consider the Dehn filling $M$ of $M_\varphi$ along the slopes $\gamma_{lat}$ in each boundary torus. Recall that this $(S,\varphi)$ is called an abstract open book decomposition for $M$. There are results relating properties of $(S,\varphi)$ and $M$. For example, the following result is an analog of the Nielsen-Thurston classification for open book decompositions.

\begin{theorem}\label{theorem:20} \cite[Theorem 8.3]{IK17} Let $(S,\varphi)$ be an abstract open book decomposition of $M$ with $|\mathcal{fD}(\varphi,R)|> 4$ for every boundary component $R \subset \partial S$. Then:
\begin{enumerate}
\item $M$ is toroidal if and only if $\varphi$ is (isotopic to a) reducible map,
\item $M$ is hyperbolic if and only if $\varphi$ is (isotopic to a) pseudo-Anosov map,
\item $M$ is Seifert fibred if and only if $\varphi$ is (isotopic to a) periodic map.
\end{enumerate}
\end{theorem}

\noindent
Thus theorem \ref{theorem:20} (2) gives conditions under which $M$ is hyperbolic if and only if the mapping torus $M_\varphi$ is hyperbolic. Theorem \ref{theorem:3} gives a similar result with additional control on the geometry of $M$. We may rephrase theorem \ref{theorem:3} as follows:\\

\noindent
{\bf Theorem \ref{theorem:3}.} Fix $0 < \varepsilon \leq \log(3), J > 1$ and let $(S, \varphi)$ be an abstract open book decomposition of a 3-manifold $M$ such that:
\begin{itemize}
\item[-] For every component $R \subset \partial S$ we fill along: \\ \indent\indent\indent\indent\indent\indent $|\mathcal{fD}(\varphi,R)| > F(S,\varepsilon,J) := K(\varepsilon,J)^2\cdot (2+|\chi(S)|) \cdot(1+19296\cdot\chi(S)^6)+3$;
\item[-] $\varphi$ is isotopic to a pseudo-Anosov map.
\end{itemize}
Then $M$ is hyperbolic and there are $J$-bi-Lipschitz inclusions $M^{\geq \varepsilon}\inject M_{\varphi}^{\geq \varepsilon/1.2}, M_\varphi^{\geq \varepsilon}\inject M^{\geq \varepsilon/1.2}$. (Here $M^{\geq \varepsilon}$ denotes the $\varepsilon$-thick part of $M$.)

\begin{remark}\label{remark:21} In theorem \ref{theorem:3} we may take $K(\varepsilon,J)$ to be the function given in \cite[Corollary 9.34]{FPS22}, namely:
$$K(\varepsilon, J)^2=\max\left\{ \frac{2\pi \cdot 6771\cosh^5(0.6\varepsilon + 0.1475)}{\varepsilon^5}+11.7, \frac{2\pi \cdot 11.35}{\varepsilon^{5/2}\log(J)} + 11.7\right\}.$$
\end{remark}

\noindent
{\it Proof:} For each boundary component $R \subset \partial S$ let $C_R$ denote the corresponding maximal cusp in $M_\varphi$ and let $\gamma_{lat,R}, \gamma_{\partial ,R}, \mathcal{sk}(\varphi, R),$ and $\text{height}(\partial C_R)$ be the quantities corresponding to the geometry of $C_R$ as given in definition \ref{definition:1}. Let the binding of $(S, \varphi)$ have $b$ components. In particular $b\leq |\chi(S)|+2$. \\

\noindent
It is shown in \cite[Corollary 9.34]{FPS22} that if for each boundary component $R \subset \partial S$ that we fill along we have $\hat{L}(\gamma_{lat,R}) \geq \sqrt{b}\cdot K(\varepsilon,J)$, then $M$ is hyperbolic and the desired $J$-bi-Lipschitz inclusions exist. For a component $R \subset \partial S$ we fill along, let us show that $\hat{L}(\gamma_{lat,R}) \geq \sqrt{b}\cdot K(\varepsilon,J)$. For this observe that:

$$\hat{L}(\gamma_{lat,R}) \geq \frac{\text{height}(\partial C_R)}{\sqrt{\text{height}(\partial C_R) \cdot l(\gamma_{\partial ,R})}}=\frac{\sqrt{\text{height}(\partial C_R)}}{\sqrt{l(\gamma_{\partial , R})}} \geq \frac{\sqrt{d_{\mathcal{A}(S,R)}(\varphi)}}{\sqrt{536\chi(S)^4 \cdot l(\gamma_{\partial ,R}),
}},$$

\noindent
where the second inequality is due to theorem \ref{theorem:10}. Hence we may, without loss of generality, assume that $\text{height}(\partial C_R) \leq b\cdot K(\varepsilon, J)^2 \cdot l(\gamma_{\partial ,R})$ and $d_{\mathcal{A}(S,R)}(\varphi)\leq b\cdot K(\varepsilon, J)^2\cdot 536\cdot \chi(S)^4\cdot l(\gamma_{\partial, R})$. In particular by lemma \ref{lemma:9}, $d_{\mathcal{A}(S,R)} \leq b\cdot K(\varepsilon, J)^2\cdot 3216\cdot |\chi(S)|^5$. Moreover, we have:
$$\hat{L}(\gamma_{lat,R}) =\frac{l(\gamma_{lat,R})}{\sqrt{\text{height}(\partial C_R) \cdot l(\gamma_{\partial ,R})}}\geq \frac{|\mathcal{sk}(\varphi, R)|\cdot l(\gamma_{\partial ,R})}{\sqrt{\text{height}(\partial C_R) \cdot l(\gamma_{\partial ,R})}} \geq \frac{|\mathcal{sk}(\varphi, R)|}{\sqrt{b}\cdot K(\varepsilon, J)}.$$
\noindent
Again, we may thus assume $|\mathcal{sk}(\varphi,R)| \leq b\cdot K(\varepsilon,J)^2$. By theorem \ref{theorem:2}, we know that $|\mathcal{fD}(\varphi,R)| \leq |\mathcal{sk}(\varphi,R)|+6|\chi(S)|\cdot d_{\mathcal{A}(S,R)}(\varphi, R) + 3$. Combining this with the previous bounds we see $|\mathcal{fD}(\varphi, R)|\leq F(S, \varepsilon, J)$. $\square$\\

\noindent
Theorem 3, in particular, allows us to related the volume of $M_\varphi$ and its Dehn fillings. For completeness, we now record such a result:

\begin{proposition}\label{proposition:22} Consider a surface $S$, a constant $c > 2\pi$, and a map $\varphi \in \text{Diff}^+(S, \partial S)$ isotopic to a pseudo-Anosov map. Let $M_\varphi$ denote the corresponding mapping torus and for each component $R \subset \partial S$ let
 $C_R \subset M_\varphi$ be the corresponding maximal cusp with associated quantities $\gamma_{lat,R},\gamma_{\partial,R},\mathcal{sk}(\varphi,R),$ and $\text{height}(\partial C_R)$ as described in definition \ref{definition:1}. Let M be the Dehn filling of $M_\varphi$ along some of the slopes 
$\gamma_{lat,R}$ for $R \subset \partial S$ and suppose in each component $R$ we fill along: $|\mathcal{fD}(\varphi, R)| > c \cdot (1 + 3216 \cdot |\chi(S)|^5)+3.$ Then, we have:
$$\text{vol}(M_\varphi) \geq \text{vol}(M) \geq (1 - (2\pi/c)^2)^{3/2} \cdot \text{vol}(M_\varphi).$$
\end{proposition}

\noindent
{\it Proof:} The upper bound is immediate from Thurston's hyperbolic Dehn surgery theorem. The lower bound would follow from \cite[Theorem 1.1]{FKP08} if in each boundary component $R$ we fill along $l(\gamma_{lat,R}) \geq c$. Fix some $R\subset \partial S$ and show this bound holds for $R$. Observe that $l(\gamma_{lat,R}) \geq |\mathcal{sk}(\varphi, R)|$. So we may assume $|\mathcal{sk}(\varphi, R)| \leq c$. Similarly, by theorem \ref{theorem:10}, $l(\gamma_{lat,R}) \geq \text{height}(\partial C_R) \geq d_{\mathcal{A}(S,R)}(\varphi)/(536 \cdot \chi(S)^4)$. 
So we may assume $d_{\mathcal{A}(S,R)}(\varphi) \leq c \cdot 536 \cdot \chi(S)^4$. By these assumptions and theorem \ref{theorem:2}, for each $R \subset \partial S$ that we fill along:
$$|\mathcal{fD}(\varphi,R)|\leq|\mathcal{sk}(\varphi,R)|+6\cdot|\chi(S)|\cdot d_{\mathcal{A}(S,R)}(\varphi)+3<c+6|\chi(S)|\cdot c\cdot 536 \cdot \chi(S)^4+3.$$
This completes the proof. $\square$ \\

\noindent
Finally, braids complements are a class of 3-manifolds which lend themselves well to the above techniques. The {\it braid group} $B_n$ is given by:
$$B_n := \langle \sigma_1, . . . , \sigma_{n-1}|\sigma_i\sigma_k = \sigma_k\sigma_i\text{ for }|i - k| \geq 2, \sigma_i\sigma_{i+1}\sigma_i = \sigma_{i+1}\sigma_i\sigma_{i+1}\rangle.$$
\noindent
We may interpret elements of $B_n$ as braids on $n$-strands as shown in figure \ref{figure:15} left, such that $\sigma_i$ corresponds to the $i$-th strand crossing under the $(i + 1)$-st strand, and group multiplication corresponding to concatenation of braids. Hence $B_n$ is the mapping class group relative boundary of the $n$-punctured disc $D_n$. Every element $\beta \in B_n$ gives rise to a knot $\hat{\beta} \subset S^3$ called the {\it braid closure} as shown in figure \ref{figure:15} right.\\

{\centering
\includegraphics[width=7.5cm]{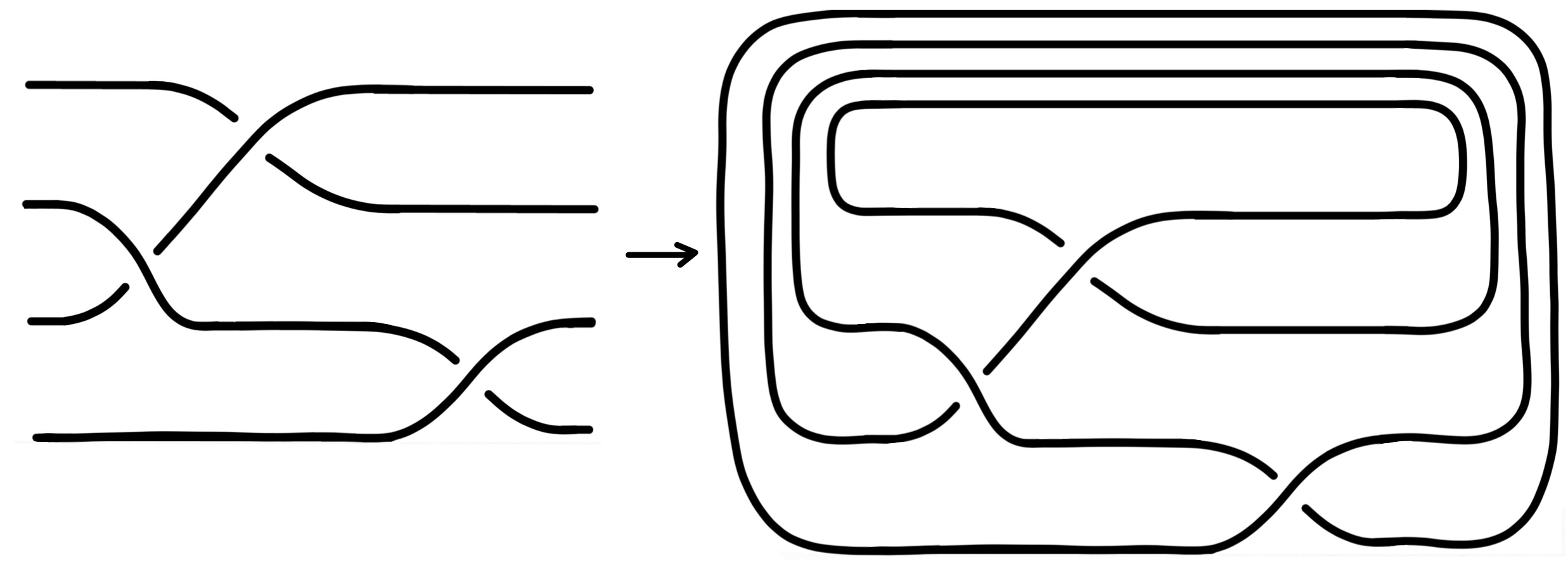}
\captionof{figure}{The braid $\sigma_1^{-1}.\sigma_3^{-1}.\sigma_2$ and its closure}
\label{figure:15}}
\vskip 0.3cm

\noindent
Given $\beta \in B_n$ we may associate with it a real number $\lfloor \beta\rfloor_D$ called the {\it Dehornoy floor}. The Dehornoy floor has a purely combinatorial definition and has been extensively studied, see for example \cite{DDRW08}. Moreover, interpreting $\beta \in B_n$ as an element of $\text{MCG}^+(D_n,S^1)$ we may define $\mathcal{fD}(\beta) = \mathcal{fD}(\beta,S^1)$, for $S^1$ the boundary of $D_n$. Then, by \cite[Lemma 7.2.1]{Mal05}, $\lfloor \beta\rfloor_D \leq \mathcal{fD}(\beta) \leq \lfloor \beta\rfloor_D + 1$ and, by \cite[Theorem 7.3]{Mal05}, $\mathcal{fD}(\beta) = \lim_{k\to\infty}\lfloor \beta^k\rfloor_D/k$. Combining this with theorem \ref{theorem:3} allows us to study the hyperbolic geometry of braid closures.\\

\noindent
Let $\mathcal{T}=\{(\sigma_i\cdots \sigma_{j-1} )(\sigma_i\cdots \sigma_{j-2} )\cdot (\sigma_i\cdot \sigma_{i+1} )\cdot \sigma_i|1\leq i<j\leq n\}$ and $l_\mathcal{T}(\rho):=\min\{l|\rho=\rho_1^{a_1}\cdots \rho_l^{a_l},\rho_i \in\mathcal{T}\}$. Note that $\mathcal{T}$ is a generating set for $B_n$ with elements as shown in figure \ref{figure:16}.\\

{\centering
\includegraphics[width=5cm]{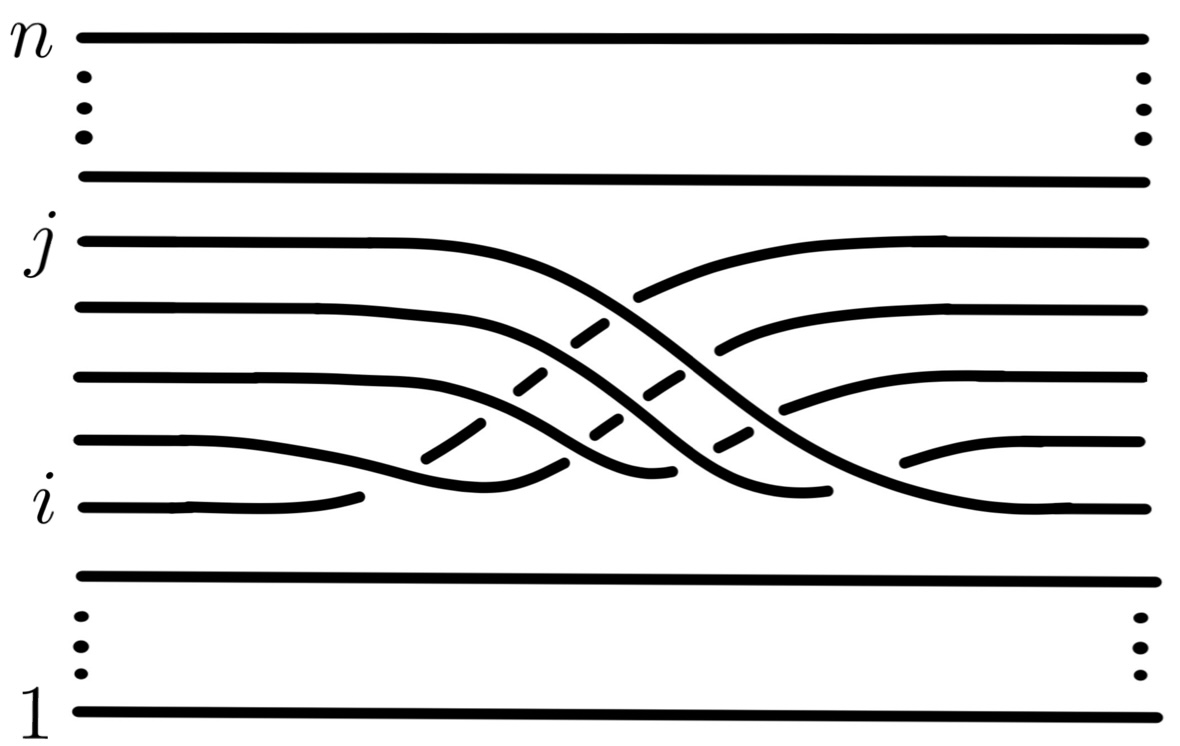}
\captionof{figure}{The braid $(\sigma_i\cdots \sigma_{j-1} )(\sigma_i\cdots \sigma_{j-2} )\cdot (\sigma_i\cdot \sigma_{i+1} )\cdot \sigma_i\in B_n$.}
\label{figure:16}}
\vskip 0.3cm

\noindent
{\bf Theorem 4.} There exists $A = A(n) \in R$ such that for a braid $\beta \in B_n$ with hyperbolic braid closure $\hat{\beta}$ and $|\lfloor \beta\rfloor_D| \geq 7 \cdot (1 + 3216 \cdot (n - 1)^5)+4$ we have:
$$\frac 1A \cdot \min_{\rho\in \text{Conj}(\beta)}l_\mathcal{T}(\rho) \leq \text{vol}(S^3 - \hat{\beta}) \leq A \cdot \min_{\rho\in\text{Conj}(\beta)}l_\mathcal{T}(\rho),$$
for $\text{Conj}(\beta) \subset B_n$ the conjugacy class of $\beta$.\\

\noindent
{\it Proof:} In an abuse of notation, we may view $\hat{\beta}$ as a link in the solid torus $D^2 \times S^1 \subset S^3$. Then $D^2 \times S^1 - \hat{\beta}$ is the mapping torus of $\beta \in \text{MCG}^+(D_n, S^1)$. By theorem \ref{theorem:20}, we know $\beta$ is isotopic to a pseudo-Anosov map. It is shown in \cite[Proposition 3.2.10]{Cap16} that there exists $A' = A'(n)$ such that:
$$\frac 1{A'} \cdot \min_{\rho\in\text{Conj}(\beta)}l_\mathcal{T}(\rho) \leq \text{vol}(D^2 \times S^1 - \hat{\beta}) \leq A' \cdot \min_{\rho\in\text{Conj}(\beta)}l_\mathcal{T}(\rho).$$
The braid closure complement $S^3 - \hat{\beta}$ is a Dehn filling of $D^2 \times S^1 - \hat{\beta}$ along $\gamma_{lat, S^1}$. Hence, the result follows from proposition \ref{proposition:22} with $c = 7$. $\square$

\section{Open Book Decompositions of a Fixed Manifold} \label{section:8}

\noindent
The following theorem of M. Hedden and T. Mark was proved using the machinery of Heegaard Floer homology.

\begin{theorem}\label{theorem:23}
 Given a closed oriented manifold $M$, there exists $N_M \in  \mathbb N$ such that any abstract open book
decomposition $(S, \varphi)$ of $M$ with connected binding has $|\mathcal{fD}(\varphi)| \leq  N_M + 1$.
\end{theorem}

\noindent
We will use theorem \ref{theorem:2} to prove a similar result. Theorem \ref{theorem:5} is theorem \ref{theorem:23} with $(S, \varphi)$ now required to have bounded $\chi(S)$, but no longer required to have connected binding. We will make use of JSJ-decompositions, so let us briefly record the following definition and classical facts:\footnote{Our definition of JSJ-decomposition is stricter than the classical definition, as we require pieces to be hyperbolic or SFS rather than simple or SFS. In light of the Perelman-Thurston Geometrization conjecture these definitions are equivalent.}

\begin{definition} \cite[Proposition V.5.2]{JS79} A collection of disjoint incompressible tori in a 3-manifold $M$ is a {\it JSJ-decomposition} if:
\begin{itemize}
\item[-] it cuts $M$ into hyperbolic or Seifert fibred components, none of which is a solid torus;
\item[-] no two adjacent Seifert fibred components admit Seifert fibrations with isotopic fibers;
\item[-] if any component is a thickened torus $T^2 \times  [0, 1]$, then its two boundary components are the same JSJ-torus, and the homeomorphism $T^2\to T^2$ induced by pushing the gluing map through the product structure is not periodic.
\end{itemize}
Equivalently a JSJ-decomposition is a minimal collection of incompressible tori cutting M into hyperbolic or Seifert fibred components. We say a manifold $M$ {\it admits} a JSJ-decomposition if $M$ contains a collection of incompressible tori cutting $M$ into Seifert fibred and hyperbolic pieces. Prime manifolds are known to admit a JSJ-decomposition. Moreover, when a JSJ-decomposition exists, it is known to be unique up to isotopy. \end{definition}

\noindent
For an introduction to Seifert fibred spaces see \cite{Mar22}. We recall that a Seifert fibration, up to isomorphism, can be expressed via the data $[S,q_1/p_1,...,q_k/p_k]$, where $S$ is a surface and $p_i\neq 0$ are the multiplicities of the corresponding fibres (unlike some authors we allow $p_i=1$). We record the following classical facts about Seifert fibrations:

\begin{theorem}\label{theorem:25} \cite[Theorem 10.1]{Wal67} \cite[Lemma 3]{Tol78} Seifert fibred spaces admit unique Seifert fibrations up to isomorphism, aside from the following:
\begin{itemize}
\item[-] the solid torus fibres as $[D^2]$ and $[D^2, q/p]$;
\item[-] the oriented circle bundle over the M\"obius band also fibres as $[D^2, 1/2, -1/2]$;
\item[-] $[S^2, 1/2, -1/2, q/p]$ and $[\mathbb RP^2, p/q]$ are Seifert fibrations for the same space;
\item[-] the oriented circle bundle over the Klein bottle also fibres as $[S^2, 1/2, 1/2, -1/2, -1/2]$;
\item[-] lens spaces (including $S^2 \times  S^1$) fibre in many ways.
\end{itemize}
Moreover, suppose M is an irreducible, $\partial $-irreducible Haken manifold with boundary, which is not $T^3, T^2 \times  [0, 1]$, nor a twisted circle or interval bundle over the Klein bottle. Then any two Seifert fibrations of $M$ are related by an isomorphism isotopic to the identity. 
\end{theorem}

\begin{proposition}\label{proposition:26} \cite[Corollary 10.4.10]{Mar22} Every Seifert fibred space is irreducible and $\partial $-irreducible, except $S^2 \times  S^1, D^2 \times  S^1$, and $\mathbb RP^2\tilde{\times} S^1$.
\end{proposition}

\noindent
{\bf Theorem \ref{theorem:5}.} Given an oriented manifold $M$, that is not a lens space or a solid torus, there exists $N_{M,\chi}$ such that any abstract open book decomposition $(S, \varphi)$ of $M$ with $\chi(S) = \chi$ has $|\mathcal{fD}(\varphi, R)| \leq  N_{M,\chi} + 1$ for some component $R \subset  \partial S$.\\

\noindent
{\it Proof:} If M admits a JSJ-decomposition, let $\text{sys}(M)$ denote the length of the shortest geodesic in a hyperbolic JSJ-component of $M$, and let $\text{sing}(M)$ denote the maximal multiplicity of a singular fibre in any Seifert fibration of a JSJ-component of $M$. This is well-defined by theorem \ref{theorem:25}, since no JSJ-component of $M$ is a solid torus or lens space. Call a manifold {\it bad} if it is hyperbolic with a geodesic of length less than $\text{sys}(M)$ or has a Seifert fibration with singular fibre of multiplicity at least $\max\{3,\text{sing}(M)\}$, and is not a lens space or solid torus. By construction, $M$ does not admit a JSJ-decomposition with a bad component. Moreover, by theorem \ref{theorem:25}, a bad Seifert fibred manifold with boundary has a unique Seifert fibration up to isotopy.\\

\noindent
Fix an abstract open book decomposition $(S, \varphi)$ of $M$ with $\chi(S) = \chi$. Consider a, possibly empty, essential multi-curve $\Lambda  \subset  S$ setwise preserved by $\varphi$ such that $S\setminus \Lambda $ is a disjoint union $S_1 \sqcup \ldots \sqcup  S_n$ of, possibly disconnected, surfaces such that on each $S_i$ the first return map is pseudo-Anosov or periodic, but not reducible. Let us denote $\varphi_i := \varphi|_{S_i}$. The mapping torus $M_\varphi = M\setminus B$ can be cut along the tori formed by the suspension of $\Lambda $ into pieces $M_{\varphi_1}, \ldots,  M_{\varphi_n}$ . Those mapping tori $M_{\varphi_i}$ with $\varphi_i$ pseudo-Anosov will be hyperbolic and those with $\varphi_i$ periodic will be Seifert fibred. Note that the Dehn filling of $M_\varphi$ to $M$ restricts
to Dehn fillings of the $M_{\varphi_i}$ to manifolds that we denote by $\widehat{M_{\varphi_i}}$.\\

\noindent
Suppose $\varphi_i$ is pseudo-Anosov and we Dehn fill along some boundary torus $C_R$ corresponding to $R\subset \partial S_i$. Let $\gamma_{lat,R}, \gamma_{\partial ,R},$ and $\text{height}(\partial C_R)$ be the quantities describing the geometry of $\partial C_R$ as given in definition \ref{definition:1}. By \cite[Corollary 6.13]{FPS22}, there is a $K(\text{sys}(M)) \geq  0$ such that after Dehn filling any hyperbolic manifold along slopes with normalised length at least $K(\text{sys}(M))$ the manifold remains hyperbolic and has a geodesic shorter than $\text{sys}(M)$. Observe that the normalised length $\hat{L}(\gamma_{lat,R})$ satisfies:

$$\hat{L}(\gamma_{lat,R})=\frac{l(\gamma_{lat,R} )}{\sqrt{\text{height}(\partial C_R) \cdot  l(\gamma_{\partial ,R})}}\geq \sqrt{\frac{l(\gamma_{lat,R})}{l(\gamma_{\partial, R})}}\cdot \sqrt{\frac{l(\gamma_{lat,R})}{\text{height}(\partial C_R)}}\geq \sqrt{|\mathcal{sk}(\varphi_i, R)|}.$$

\noindent
Here the final inequality is due to $l(\gamma_{lat,R}) \geq \text{height}(\partial C_R)$. Similarly, we observe that:

$$\hat{L}(\gamma_{lat,R})\geq\frac{\text{height}(\partial C_R)}{\sqrt{\text{height}(\partial C_R) \cdot  l(\gamma_{\partial ,R})}}\geq \frac{\sqrt{d_{\mathcal{A}(S_i,R)}(\varphi_i)}}{\sqrt{536 \cdot  \chi(S_i)^4 \cdot  6 \cdot  |\chi(S_i)|}},$$

\noindent
where the final inequality is due to theorem \ref{theorem:10} and lemma \ref{lemma:9}. Finally, recall that by theorem \ref{theorem:2}, $|\mathcal{fD}(\varphi_i, R) - \mathcal{sk}(\varphi_i, R)| \leq  6 \cdot  |\chi(S_i)| \cdot  d_{\mathcal{A}(S_i,R)}(\varphi_i) + 3.$ In summary, if for each $R \subset  \partial S_i$ that we fill along we have:
$$|\mathcal{fD}(\varphi,R)|\geq K(\text{sys}(M))^2 +3+6\cdot |\chi(S)|\cdot K(\text{sys}(M))^2 \cdot 536\cdot \chi(S)^4 \cdot 6\cdot |\chi(S)|,$$
then $\widehat{M_{\varphi_i}}$ will be hyperbolic with a geodesic of length less than $\text{sys}(M)$.\\

\noindent
Suppose $\varphi^P_i$ is isotopic, not relative boundary, to $\text{id}: S_i \to S_i$ for some $P \in  \mathbb N$ and we Dehn fill along some boundary torus $C_R$ corresponding to $R\subset\partial S_i$. Then the $P$-fold cyclic cover of $M_{\varphi_i}$ is $S_i \times  S^1$ and this circle-bundle structure descends to a Seifert fibration of $M_{\varphi_i}$. Let $\gamma_{fib,R} = \gamma_{lat,R} + D \cdot  \gamma_{\partial ,R} \in  \pi_1(\partial C_R)$ be a fibre of this Seifert fibration. Then $\varphi^P_i$ must be Dehn twisted $D$ times in $R$ along the isotopy to $\text{id}$ and hence $|\mathcal{fD}(\varphi_i,R)| = |D|/P$. Hence, if in each component $R$ that we fill along: $\mathcal{fD}(\varphi_i,R)\neq 0$, then $\widehat{M_{\varphi_i}}$ is Seifert fibred with the core curves of the Dehn filling corresponding to singular fibres of multiplicity $|\mathcal{fD}(\varphi_i, R)| \cdot  P$.\\

\noindent
Hence we may take $N_{M,\chi}$ large enough that if in each boundary component $R \subset  \partial S$ we fill along: $|\mathcal{fD}(\varphi, R)| > N_{M,\chi} + 1$, then some $\widehat{M_{\varphi_i}}$ is bad. Let us show that this contradicts our definition of $\text{sys}(M)$ and $\text{sing}(M)$.\\

\noindent
Since $\Lambda $ is essential, one of the $\widehat{M_{\varphi_i}}$ can only be a solid torus if $M$ itself is a solid torus, a contradiction. Thus by proposition \ref{proposition:26} the boundary tori of the $\widehat{M_{\varphi_i}}$ are incompressible within $M$. Moreover, none of the $\widehat{M_{\varphi_i}}$ are a thickened torus $T^2 \times  [0, 1]$. If no two Seifert fibred $\widehat{M_{\varphi_i}}$ have Seifert fibrations with fibred slopes lining up, then $\widehat{M_{\varphi_i}}$ give a JSJ-decomposition of $M$ and the contradiction is immediate.\\

\noindent
Suppose two Seifert fibred $\widehat{M_{\varphi_i}}$ have Seifert fibrations with fibres lining up. Denote these $A$ and $B$. If one of $A, B$ is bad, then it has a unique Seifert fibration up to isotopy and thus the union of Seifert fibrations on $A, B$ gives a bad Seifert fibration of $A \cup B$. Thus replacing $A, B$ with $A \cup B$ again gives a decomposition of $M$ into hyperbolic and Seifert fibred components, at least one of which is bad. By induction on the number of JSJ-components, we may again assume no two Seifert fibred components have fibers lining up and are done. $\square$\\

\noindent
This theorem is not true for $M$ a lens space or a solid torus, as these have Seifert fibrations with singular fibres of arbitrarily high multiplicity. Finally, one may hope to drop the assumption that $\chi(S) = \chi$ is fixed. To do so, one would need to remove the $|\chi(S)|$ term from theorem \ref{theorem:2} and find a $\chi(S)$ independent bound on $d_{\mathcal{A}(S,R)}$. The former would follow if one could replace the upper bound of lemma \ref{lemma:9} by a multiple of $l(\gamma_\partial )$. The latter (for connected binding) is a conjecture of S. Schleimer and has been shown, among others, to hold for Haken manifolds \cite{Cen22}. We record this as a question.

\begin{question} Given an oriented manifold $M$, that is not a lens space or solid torus, is there a constant $N_M$ such that for any abstract open book decomposition $(S, \varphi)$ of $M$ there is a boundary component $R \subset  \partial S$ such that $|\mathcal{fD}(\varphi, R)| \leq  N_M$?\\
\end{question}

\section{Applications to Contact Structures}\label{section:9}

\noindent
We may use fractional Dehn twist coeffients and theorem \ref{theorem:2} to relate the contact topology and the geometry of hyperbolic 3-manifolds. Recall theorem \ref{theorem:6}, which gives a geometric criteria on the binding of an open book decomposition $(B, \pi)$, for the corresponding contact structure on $M$ to be tight.\\

\noindent
Consider an oriented geodesic $\gamma$ in an oriented hyperbolic manifold $M$, take a normal vector $n \in  T_pM$ to $\gamma$ at $p \in  \gamma$, and parallel transport $n$ along $\gamma$ to obtain another normal vector $n' \in  T_pM$. The {\it torsion} $\tau(\gamma) \in  (-\pi,\pi]$ is the counter-clockwise angle from $n$ to $n'$. We may rephrase theorem \ref{theorem:6} as follows:\\

\noindent
{\bf Theorem \ref{theorem:6}.} Let $(\gamma, \pi)$ be an open book decomposition for a closed hyperbolic manifold $M$. Suppose $\gamma$ is a geodesic of length $l$ and torsion $\tau$. There exists positive constants $L_{max} = L_{max}(S), tor_{min} = tor_{min}(S)$ depending only on $\chi(S)$ for $S$ a page of $(\gamma, \pi)$, such that if $l^2+\tau^2 \leq L_{max}$ and $\tau/l\geq tor_{min}$, then the contact structure corresponding to $(\gamma,\pi)$ is tight.\\

\noindent
To prove theorem \ref{theorem:6}, we combine theorem \ref{theorem:2} with theorem \ref{theorem:12} and we relate the geometry of $\gamma$ to the cusp geometry of $M\setminus \gamma$. The later can be done via the angle deformation theory of Hodgeson-Kerckhoff. We now state the relevant result and prove it in section \ref{section:10}.

\begin{lemma}\label{lemma:28}
Let $M$ be a closed hyperbolic manifold with open book decomposition $(\gamma, \pi)$ with corresponding abstract open book decomposition $(S,\varphi)$ and $\gamma$ a geodesic of length $l$ and torsion $\tau$. There exists positive constants $A=A(S),B=B(S)$, depending only on $\chi(S)$, such that if $l^2+\tau^2 \leq A,|\tau|/l\geq B$, then:
$$\left|\mathcal{sk}(\varphi)-\frac{2\pi\cdot \tau}{l^2 +\tau^2}\right|, \left|\frac{\text{height}(\partial C)}{l(\gamma_\partial)}-\frac{2\pi\cdot l}{l^2 +\tau^2}\right|\leq 92,$$
here $C$ is the cusp of $M -\gamma$ and $\mathcal{sk}(\varphi),\text{height}(\partial C)$, and $l(\gamma_\partial )$ are the quantities as given in definition \ref{definition:1}.
\end{lemma}

\noindent
Let us conclude the section with a proof of theorem \ref{theorem:6}.\\

\noindent
{\it Proof of theorem \ref{theorem:6}:} Let $(S,\varphi)$ be the abstract open book decomposition corresponding to $(\gamma,\pi)$. By theorem \ref{theorem:12}, it is sufficient to show that $\mathcal{fD}(\varphi) \geq  1$. By theorems \ref{theorem:2} and \ref{theorem:10} and since $l(\gamma_\partial ) \leq  6|\chi(S)|$, we have:
\begin{eqnarray*}
\mathcal{fD}(\varphi) &\geq & \mathcal{sk}(\varphi) - 6|\chi(S)| \cdot  d_{\mathcal{A}(S,R)}(\varphi) - 3\\
& \geq & \mathcal{sk}(\varphi)-6\cdot 536|\chi(S)|^5\cdot \text{height}(\partial C)-3\\
& \geq & \mathcal{sk}(\varphi) - 36 \cdot  536|\chi(S)|^6 \cdot  \frac{\text{height}(\partial C)}{l(\gamma_\partial)} - 3.
\end{eqnarray*}

\noindent
Now take $L_{max} \leq  A, tor_{min} \geq  B$ with $A, B$ from lemma \ref{lemma:28}, so that:
\begin{eqnarray*}
\mathcal{fD}(\varphi) & \geq & \frac{2\pi\cdot \tau}{l^2+\tau^2}-36\cdot 536|\chi(S)|^6\cdot \left(\frac{2\pi \cdot l}{l^2+\tau^2}+92\right)-3-92\\
& \geq & \frac{2\pi\cdot \tau}{l^2+\tau^2}\left( 1-19296|\chi(S)|^6\cdot \frac l\tau\right)-92 \cdot (1+19296|\chi(S)|^6)-3.
\end{eqnarray*}

\noindent
Hence for $tor_{min}$ and $\tau/(l^2 + \tau^2)$ large enough, we have $\mathcal{fD}(\varphi) \geq  1$.\\

\noindent
Observe that the region $\{(l, \tau)|\tau \geq  C \cdot  (l^2 + \tau^2)\} \subset  \mathbb R^2$, see figure \ref{figure:17}, is a disc of radius $1/2C$ with centre $(0,1/2C)$. As we see in figure \ref{figure:17}, for $tor_{min}$ large enough and $L_{max}$ small enough, points $l,\tau$ with $l \geq  0, \tau/l \geq  tor_{min}, l^2 + \tau^2 \leq  L_{max}$ lie inside this disc. Hence, for $tor_{min}$ large enough and $L_{max}$ small enough, we have that $\mathcal{fD}(\varphi) \geq  1$. This concludes the proof. $\square$\\

{\centering
\includegraphics[width=4.5cm]{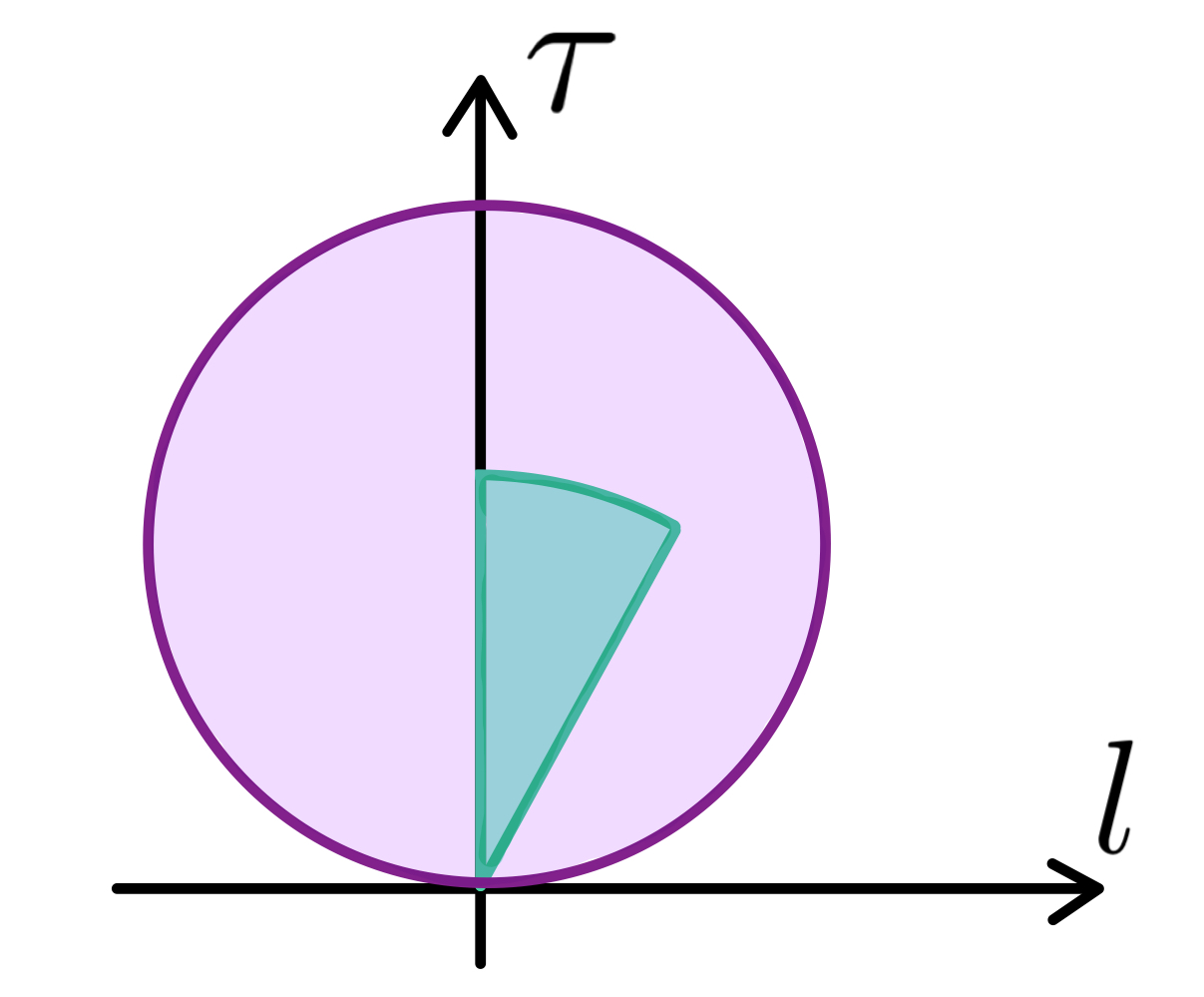}
\captionof{figure}{Domains $\{(l, \tau)|\tau \geq C \cdot (l^2 + \tau^2)\}$ in purple and $\{(l, \tau)|l > 0, \tau/l\geq tor_{min}, l^2 + \tau^2 \leq L_{max}\}$ in green for suitable $C, tor_{min}, L_{max}$.}
\label{figure:17}}
\vskip 0.3cm

\section{Angle Deformation Theory} \label{section:10}
\noindent
Consider a hyperbolic 3-manifold $M$ and a simple closed geodesic $\gamma \subset  M$. It is known \cite[Theorem 1.2.1]{Koj98}, that $M - \gamma$ admits a hyperbolic metric. In \cite{HK05} C. Hodgson and S. Kerckhoff give sufficient criteria for the complete hyperbolic structures on $M$ and $M - \gamma$ to be connected by a 1-parameter family of, so called, hyperbolic cone-manifolds structures on $M$. It is then possible to understand the change in geometry along this 1-parameter family. In this section we use these techniques to prove lemma \ref{lemma:28}. Similar work relating the cusp geometry of $M-\gamma$ and the length and torsion of $\gamma\in M$ has appeared in \cite{Magid_Kleinian} and \cite{Purcell-Thesis}. In fact, with some additional work, \cite[Section 3]{Purcell-Thesis} could be used for an alternative proof of Lemma \ref{lemma:28}. \\ 

\noindent
See \cite[Section 2]{FPS22} for a detailed introduction to cone-manifolds. Let $\hat{\mathbb H^3}$ be the metric completion of the universal cover of $\mathbb H^3 - \sigma$, for $\sigma$ some geodesic. Let $\hat{\sigma}$ be the points added in the metric completion, so that $\hat{\mathbb H^3}$ is an infinite cyclic branch cover of $\mathbb H^3$ with branch set $\hat{\sigma}$. Then $\mathbb C$, as an additive group, acts on $\hat{\mathbb H^3}$ by letting $l + i\tau$ ``translate'' $\hat{\sigma}$ by $l$ and ``rotate'' $\hat{\sigma}$ by $\tau$ . We denote such an isometry by $\varphi_{l+i\tau}$.

\begin{definition} A {\it model solid torus} $N_{\alpha,l,\tau}$ is the quotient of $\hat{\mathbb H^3}$ by $\langle\varphi_{i\alpha}, \varphi_{l+i\tau}\rangle$ for $\alpha, l > 0$. The quotient of $\hat{\sigma}$ is the {\it core geodesic} of $N_{\alpha,l,\tau}$ which we say is a {\it cone-singularity of angle} $\alpha$ and of {\it complex length} $l + i\tau$. Notice that $\tau$ is only well-defined modulo $\alpha$, and unless stated otherwise $-\frac{\alpha}2 < \tau \leq \frac{\alpha}2$. The case $\alpha = 2\pi$ corresponds to the non-singular hyperbolic case. A {\it hyperbolic cone-manifold} is a space locally isometric to a model solid-torus. Points isometric to the core geodesic of a model solid-torus with $\alpha\neq 2\pi$ are called {\it singular}.
\end{definition}

\noindent
For $R > 0$, a model tube $U_{\alpha,l,\tau}$ is the open $R$-neighbourhood of the core curve in $N_{\alpha,l,\tau}$. The metric tensor in $\hat{\mathbb H^3}$ is $dr^2 + \text{cosh}^2(r)d\xi^2 + \text{sinh}^2(r)d\theta^2$, where $r$ is the distance from $\hat{\sigma}, \xi$ is the distance along $\hat{\sigma}$, and $\theta$ is the angle around $\hat{\sigma}$. One may compute that the closed model tube $U_R := U_{\alpha,l,\tau}$ of radius $R$ has boundary a Euclidean torus. As detailed in \cite[p. 388]{HK05} this Euclidean metric can be built from a Euclidean cylinder of circumference $m$ and height $h$ by attaching ends with a twist of distance $tw$, for $m, h,$ and $tw$ as below:
$$m = \alpha \cdot \text{sinh}(R), h = l\cdot  \text{cosh}(R), tw = \tau \cdot \text{sinh}(R).$$
\noindent
Moreover, in this setup, the cylinder boundary of length $m$ corresponds to the meridian, that is the homotopically
trivial, curve of $U_R$.\\

\noindent
The following theorem summarises the requisite technical input from angle deformation theory. The bi-Lipschitz distance between two metric space $X, Y$ is the infimum of $|\log \text{lip}(f)| + | \log \text{lip}(f^{-1})|$ over all bi-Lipschitz mappings $f : X \to Y$ , where $\text{lip}$ denotes the Lipschitz constant. Define this distance to be infinite, if $X, Y$ are not bi-Lipschitz.

\begin{theorem}\label{theorem:30} \cite[Theorem 1.2, Corollary 6.3, Remark 5.8, Theorem 4.4, Equation 19, Equation 18]{HK05} Let $M$ be closed hyperbolic manifold containing a geodesic $\gamma$ of length $l < 0.111$. Then $M$ and $M - \gamma$ can be connected by a 1-parameter family $M_\alpha, \alpha \in  (0, 2\pi]$ of hyperbolic cone-manifolds homeomorphic to $M$ with cone angle $\alpha$. In $M_\alpha,\alpha \in  (0,2\pi]$ let $R(\alpha),l(\alpha),$ and $\tau(\alpha)$ denote the radius of a maximal tubular neighbourhood $U_{R(\alpha)}$ about $\gamma$, the length of $\gamma$, and the torsion of $\gamma$. Let $M_0$ denote $M - \gamma$ with the complete cusped hyperbolic metric and $U_{R(0)}$ denote its maximal cusp. The family $M_\alpha$ satisfies the following:
\begin{enumerate}
\item $M_\alpha \setminus U_{R(\alpha)}, \alpha \in  [0, 2\pi]$ is a continuous 1-parameter family with respect to the bi-Lipschitz topology;\footnote{Observe that $M_\alpha$ is a continuous family for $\alpha \in  (0, 2\pi]$ but not for $\alpha \in  [0, 2\pi]$ as $M_0$ and $M_{2\pi}$ are not even homeomorphic. Therefore, to study the limit $\alpha \to 0$, we instead consider $M_\alpha \setminus U_{R(\alpha)}$.}
\item For $\alpha \in  (0, 2\pi]$ the following holds:
$$R(\alpha)\geq 0.531, \lim_{\alpha\to 0}R(\alpha)=\infty,\text{ and }\alpha\cdot l(\alpha)\geq 3.3957\cdot \frac{\tanh(R(\alpha))}{\cosh(2R(\alpha))};$$
\item $l(\alpha), \tau (\alpha), \alpha \in  (0, 2\pi]$ are smooth and satisfy:
$$\frac{dl(\alpha)}{d\alpha}=\frac{l(\alpha)}{\alpha}(1+4\alpha^2 x)\text{ and }\frac{d\tau(\alpha)}{d\alpha}=\frac{\tau(\alpha)}{\alpha}+4\alpha y l(\alpha), \text{ with } |x|, |y|\leq \frac 1{4\alpha^2 \sinh^2(R(\alpha))}.$$
\end{enumerate}
\end{theorem}

\noindent
In the above theorem we must be careful when defining $\tau(\alpha)$ as it is initially only defined modulo $\alpha$. As we wish $\tau(\alpha)$ to be continuous we set $-\pi < \tau(2\pi) \leq  \pi$ and extend $\tau(\alpha)$ by continuity. However, as we will see in the proof of lemma \ref{lemma:28}, in our setup, $-\frac{\alpha}2 < \tau(\alpha) \leq  \frac{\alpha}2$ will continue to hold for all $\alpha$.\\

\noindent
We may now recall and prove lemma \ref{lemma:28}.\\

\noindent
{\bf Lemma \ref{lemma:28}.} Let $M$ be a closed hyperbolic manifold with open book decomposition $(\gamma, \pi)$ with corresponding abstract open book decomposition $(S,\varphi)$ and $\gamma$ a geodesic of length $l(2\pi)$ and torsion $\tau(2\pi)$. There exists $A=A(S),B=B(S)\in \mathbb R_{>0}$, depending only on $\chi(S)$, such that if $l(2\pi)^2+\tau(2\pi)^2 \leq A,|\tau(2\pi)|/l(2\pi)\geq B$, then:
$$\left|\mathcal{sk}(\varphi)-\frac{2\pi\cdot \tau(2\pi)}{l(2\pi)^2 +\tau(2\pi)^2}\right|, \left|\frac{\text{height}(\partial C)}{l(\gamma_\partial)}-\frac{2\pi\cdot l(2\pi)}{l(2\pi)^2 +\tau(2\pi)^2}\right|\leq 92.$$

\noindent
{\it Proof:} Suppose $l(2\pi)^2+\tau(2\pi)^2\leq A$ and $|\tau(2\pi)|/l(2\pi)\geq B$, for $A, B$ to be determined. Take $A < (0.111)^2$, so that by theorem \ref{theorem:30} we have a 1-parameter family $M_\alpha, \alpha \in  [0, 2\pi]$ of hyperbolic cone manifolds connecting $M$ to $M -\gamma$. Let $R(\alpha),l(\alpha),$ and $\tau(\alpha)$ be as in theorem \ref{theorem:30}. Consider the quantity:
$$\omega(\alpha):=\frac{\alpha\cdot \tau(\alpha)}{l(\alpha)^2+\tau(\alpha)^2}.$$
\noindent
We differentiate $\omega(\alpha)$ using theorem \ref{theorem:30} (3) to obtain:
\begin{eqnarray*}
\frac{d\omega}{d\alpha} &= & \frac{d}{d\alpha}\left(\frac{\alpha\cdot \tau}{l^2+\tau^2}\right)\\
&=&\frac{\tau}{l^2+\tau^2}+\alpha\frac{d}{d\alpha}\left(\frac{\tau}{\tau^2+l^2}\right)=\frac{\tau}{l^2+\tau^2}+\frac{\alpha}{(l^2+\tau^2)^2}\left(\frac{d\tau}{d\alpha}\cdot (l^2+\tau^2)-\tau\frac{d}{d\alpha}(l^2+\tau^2)\right)\\
& = & \frac{\tau+\alpha\frac{d\tau}{d\alpha}}{l^2+\tau^2}-\frac{2\alpha\tau}{(l^2+\tau^2)^2}\left(l\frac{dl}{d\alpha}+\tau\frac{d\tau}{d\alpha}\right)\\
& = & \frac 1{(l^2+\tau^2)}\left(\tau+\alpha\left(\frac{\tau}{\alpha}+4\alpha y l\right)\right)-\frac{2\alpha \tau}{(l^2+\tau^2)^2}\left(\frac{l^2}{\alpha}(1+4\alpha^2 x)+\tau\left(\frac{\tau}{\alpha}+4\alpha y l\right)\right)\\
& = & \frac 1{(l^2+\tau^2)^2}\left(2 l^2 \tau+2\tau^3+4\alpha^2 y l^3+4\alpha^2y l\tau^2-2\tau l^2 -8\tau l^2  \alpha^2 x-2\tau^3-8\alpha^2 \tau^2y l\right)\\
& = & \frac{4\alpha^2}{(l^2+\tau^2)^2}(l^3y-2\tau l^2x-\tau^2ly).
\end{eqnarray*}
Taking the substitution $t = t(\alpha) = \alpha^2$ and recalling the bounds on $|x|, |y|$ from theorem \ref{theorem:30} (3), we obtain:
\begin{eqnarray*}
\left|\frac{d\omega}{dt}\right| & = &\left|\frac{1}{2\alpha}\cdot\frac{d\omega}{d\alpha}\right| \leq  \frac{2\alpha l}{4\alpha^2 \sinh^2(R)}\left|\frac{l^2+2\tau l+\tau^2}{(l^2+\tau^2)^2}\right|\\
&\leq &\frac{1}{2\alpha l  \sinh^2(R)}\left|\frac{l^2(l^2+2\tau l+\tau^2)}{(l^2+\tau^2)^2}\right|\\
&\leq &\frac{1}{2\alpha l \sinh^2(R)}\cdot \frac{3+2\sqrt{2}}{4}\leq \frac{1.4572}{2\alpha l \sinh^2(R)}, 
\end{eqnarray*}
where the penultimate inequality is obtained by maximising the absolute value of $\frac{l^2(l^2+2\tau l+\tau^2)}{(l^2+\tau^2)^2}, l, \tau\in \mathbb R.$\\

\noindent
Recall the bound on $\alpha \cdot l$ from theorem \ref{theorem:30} (2):
$$\alpha \cdot l\geq 3.3957\cdot \frac{\tanh(R)}{\cosh(2R)}=3.3957\cdot \tanh(R)\frac{\cosh^2(R)-\sinh^2(R)}{\cosh^2(R)+\sinh^2(R)}=3.3957\cdot \tanh(R)\frac{1-\tanh^2(R)}{1+\tanh^2(R)}.$$
Combining this with the above, and recalling the identity $\sinh^2(R)=\frac{\tanh^2(R)}{1-\tanh^2(R)}$, we see that: 
\begin{eqnarray*}
\left|\frac{d\omega}{dt}\right| \leq  \frac{1.4572}{2\cdot 3.3957}\cdot \frac{1+\tanh^2(R)}{\tanh(R)\cdot (1-\tanh^2(R))}\cdot \frac 1{\sinh^2(R)} = \frac{1.4572}{2\cdot 3.3957}\cdot \frac{1+\tanh^2(R)}{\tanh^3(R)}.
\end{eqnarray*}
Noting that $\frac{1+z^2}{z^3}$ is monotonically decreasing and $R\geq 0.531$, it follows that: 
$$\left|\frac{d\omega}{dt}\right| \leq  \frac{1.4572}{2\cdot 3.3957}\cdot \frac{1+\tanh^2(0.531)}{\tanh^3(0.531)}\leq 2.3089.$$
Recall that along the deformation the maximal tube radius $R(\alpha)$ tends to infinity as $\alpha\to 0$, and $M_\alpha - U_{R(\alpha)}$ converges to $M_0 - U_0$. It follows that:
$$\lim_{\alpha\to 0} \omega(\alpha)=\lim_{\alpha\to 0}\frac{tw(\alpha)\cdot m(\alpha)}{tw(\alpha)^2+h(\alpha)^2\coth^2(R(\alpha))}=\frac{tw(0)\cdot m(0)}{tw(0)^2+h(0)^2}.$$
Hence, as before:
$$\left|\frac{m(0)\cdot tw(0)}{tw(0)^2+h(0)^2}-\frac{2\pi \cdot \tau(2\pi)}{l(2\pi)^2+\tau(2\pi)^2}\right|\leq \int_0^{(2\pi)^2}\left|\frac{d\omega}{dt}\right|dt\leq (2\pi)^2 \cdot 2.3089\leq 92.$$
Recall that in the cusp $C$ of $M_0$, $m(0)$ denotes the length of the filling slope. Moreover, the Euclidean geometry on $\partial C$ is constructed by taking a Euclidean cylinder of height $h(0)$ and circumference $m(0)$ and gluing ends along a twist of length $tw(0)$. From this, we will show: $\mathcal{sk}(\varphi) = (m(0) \cdot tw(0))/(tw(0)^2 + h(0)^2)$. However, we must take care with $tw(0)$ having a modulo $m(0)$ ambiguity. Recall we chose the convention $-\pi < \tau(2\pi) \leq \pi$ and carried this through the deformation to fix a choice of $tw(0)$. We want to show $-m(0)/2 < tw(0) < m(0)/2$. Since $m(0)/tw(0) = \lim_{\alpha\to 0}\alpha/\tau$, it is sufficient to show that $|\alpha/\tau | > 2 + \varepsilon$ along the deformation for some $\varepsilon > 0$. Let's therefore embark on the following computations:
\begin{eqnarray*}
\frac{d}{d\alpha}\left(\frac{\alpha}{l}\right)&=&\frac{l-\alpha \frac{dl}{d\alpha}}{l^2}=\frac{l-l(1+4\alpha^2x)}{l^2}=\frac{-4\alpha^2x}{l} \Rightarrow \\
\left|\frac{d}{dt}\left(\frac{\alpha}{l}\right)\right|& = & \left| \frac{1}{2\alpha}\frac{d}{d\alpha}\left(\frac{\alpha}{l}\right)\right|=
\left|\frac{2\alpha x}{l}\right|\leq \left|\frac{2\alpha}{4\alpha^2 \sinh^2(R)\cdot l}\right|\leq  \frac{1}{2\cdot 3.3957}\cdot \frac{1+\tanh^2(0.531)}{\tanh^3(0.531)}\leq 1.5845 \Rightarrow \\
\frac{\alpha}{l} &\geq &\frac{2\pi}{l(2\pi)}-1.5845\cdot (2\pi)^2\geq \frac{2\pi}{\sqrt{A}}-62.5536=:C_1.
\end{eqnarray*}
We may take $A$ small enough so that $C_1 > 0$ becomes arbitrarily large. We continue with a similar estimate:
\begin{eqnarray*}
\frac{d}{d\alpha}\left(\frac{\tau}{l}\right) & = & \frac{\frac{d\tau}{d\alpha}l -\tau \frac{dl}{d\alpha}}{l^2}=\frac{\left(\frac{\tau}{\alpha}+4\alpha y l\right)l-\tau\left(\frac{l}{\alpha}+4\alpha l x\right)}{l^2}=\frac{4\alpha l(yl-x\tau)}{l^2} \Rightarrow \\
\left|\frac{d}{dt}\left(\frac{\tau}{l}\right)\right| & = & \left|\frac{2l(yl-x\tau)}{l^2}\right|\leq \frac{2l(l+|\tau|)}{l^2}\cdot\frac{1}{4\alpha^2 \sinh^2(R)}\leq \left(1+\frac{|\tau|}{l}\right)\cdot \frac{l}{\alpha}\cdot \frac{1}{2\alpha l\sinh^2(R)}\\ 
& \leq &\left(1+\frac{|\tau|}{l}\right) \cdot \frac{1}{C_1} \cdot 1.5845.
\end{eqnarray*}
\noindent
Solving the ODE: $y'(x)=\frac{-1.5845}{C_1}(1+y(x))\Rightarrow y(x)=-1+\exp\left(\frac{-1.5845x}{C_1}\right)\cdot c$, for $c$ some constant, we find: 
\begin{eqnarray*}
\left|\frac{\tau}{l}\right| & \geq & -1+\exp\left(-\frac{1.5845}{C_1}(2\pi)^2\right)\cdot \left(\frac{|\tau(2\pi)|}{l(2\pi)}+1\right).
\end{eqnarray*}
Take $A$ small enough, that $C_1$ is large enough to guarantee: $\text{exp}\left(-\frac{1.5845}{C_1}(2\pi)^2\right)\geq \frac 12$. So: 
\begin{eqnarray*}\left|\frac{\tau}{l}\right| & \geq & \frac{B-1}{2}:=C_2.\end{eqnarray*}
Finally, we compute: 
\begin{eqnarray*}
\frac{d}{d\alpha}\left(\frac{\alpha}{\tau}\right)& = &\frac{\tau-\alpha\frac{d\tau}{d\alpha}}{\tau^2}=\frac{\tau-\alpha\left(\frac{\tau}{\alpha}+4\alpha y l\right)}{\tau^2}=\frac{-4\alpha^2 y l}{\tau^2} \Rightarrow \\
\left| \frac{d}{dt}\left(\frac{\alpha}{\tau}\right)\right| & = & \left|\frac{2\alpha y l}{\tau^2}\right|\leq \left|\frac{2\alpha l}{4\alpha^2 \sinh^2(R)\tau^2}\right|=\frac{l^2}{\tau^2} \frac{1}{2\alpha l \sinh^2(R)}\leq \frac{l^2}{\tau^2}\cdot 1.5845\Rightarrow \\
\left|\frac{\alpha}{\tau}\right| & \geq & \frac{2\pi}{|\tau(2\pi)|}-\frac1{C_2^2} \cdot 62.5536.
\end{eqnarray*}
Take $A< \pi^2$ so that by assumption $|\tau(2\pi)|\leq \sqrt{A}<\pi$ and take $B$ large enough that along the deformation $|\alpha/\tau|$ remain greater than 2 and thus $-m(0)/2 < tw(0) < m(0)/2.$\\

{\centering
\includegraphics[width=4.5cm]{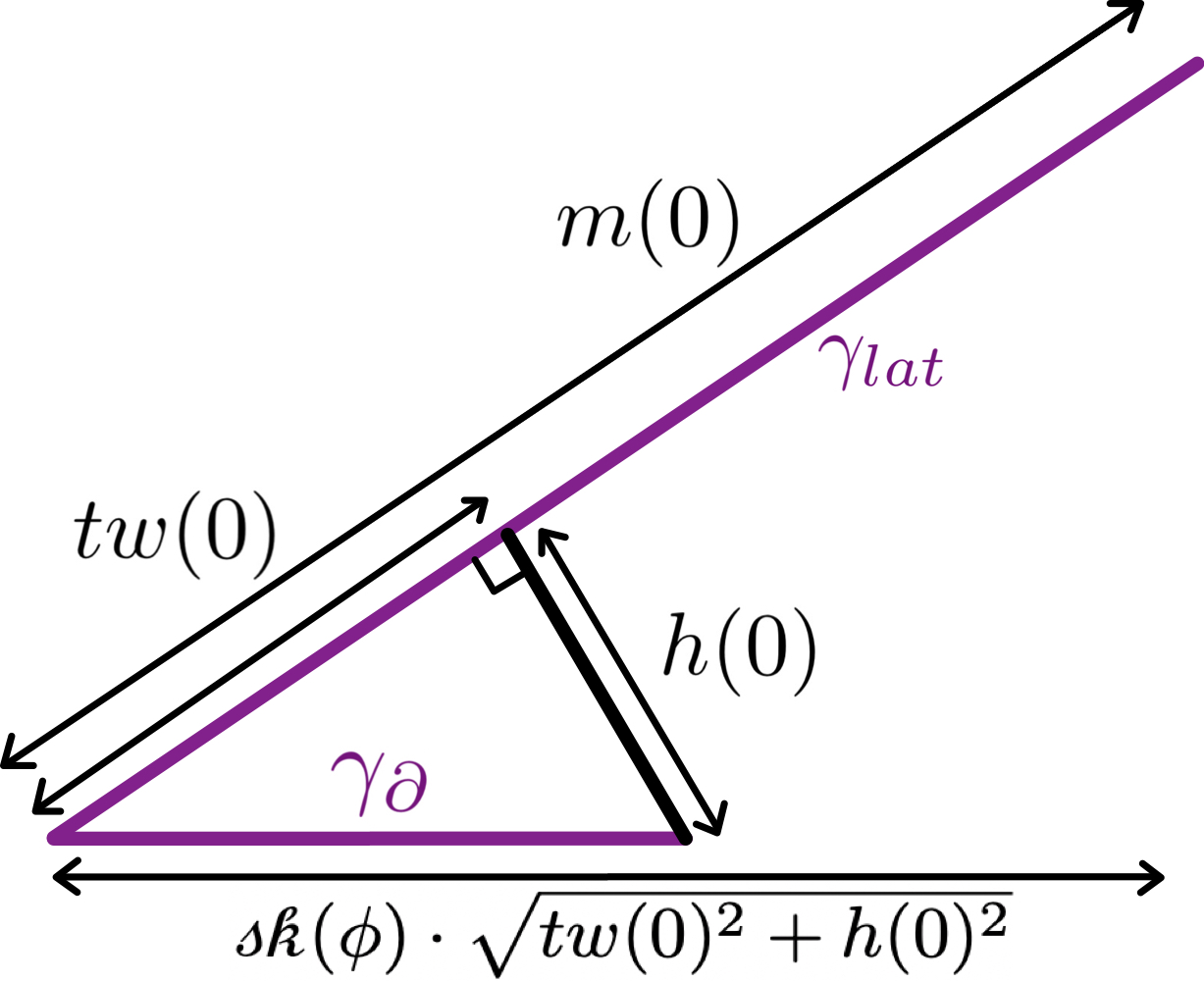}
\captionof{figure}{Curves $\gamma_\partial , \gamma_{lat}$ in $\partial C$ with segments of lengths $\mathcal{sk}(\varphi) \cdot \sqrt{tw(0)^2 + h(0)^2}, tw(0),$ and $m(0)$.}
\label{figure:18}}
\vskip 0.3cm

\noindent
Recall that $M - \gamma$ is a mapping torus $M_\varphi$. Denote, as usual, the filling slope by $\gamma_{lat}$ and the boundary of a fibre by $\gamma_\partial$. Recall in the maximal cusp $C$, we have $l(\gamma_\partial) \leq  6|\chi(S)|.$ Moreover by \cite[Corollary 6.13]{FPS22}, the normalised length $l(\gamma_{lat})/\sqrt{\text{Area}(\partial C)} \geq  \frac{2\pi}{l(2\pi)} - 16.17$ and by lemma \ref{lemma:8} $\text{Area}(\partial C) \geq \sqrt{3}/4.$ Thus for $A$ small enough in terms of $\chi(S)$, we get $l(\gamma_{lat})>2\cdot l(\gamma_\partial)$. Hence the quantities $m(0),tw(0),$ and $h(0)$ are as shown in figure \ref{figure:18}. By comparing similar triangles, it follows that:
\begin{eqnarray*}
\frac{\mathcal{sk}(\varphi)\cdot \sqrt{tw(0)^2+h(0)^2}}{m(0)}=\frac{tw(0)}{\sqrt{tw(0)^2+h(0)^2}}\Rightarrow \mathcal{sk}(\varphi)=\frac{tw(0)\cdot m(0)}{tw(0)^2+h(0)^2}.
\end{eqnarray*}
\noindent
This complete the first inequality. We proceed similarly for the second inequality.
\begin{eqnarray*}
\frac{d}{d\alpha}\left(\frac{\alpha\cdot  l}{l^2+\tau^2}\right) & = & \frac{l}{l^2+\tau^2}+\alpha \frac{d}{d\alpha}\left(\frac{l}{l^2+\tau^2}\right) =\frac{l}{l^2+\tau^2}+\frac{\alpha}{(l^2+\tau^2)^2}\left(\frac{dl}{d\alpha}(l^2+\tau^2)-l\frac{d}{d\alpha}(l^2+\tau^2)\right) \\
 & = & \frac{l+\alpha \frac{dl}{d\alpha}}{l^2+\tau^2}-\frac{2\alpha l}{(l^2+\tau^2)^2}\left(l\frac{dl}{d\alpha}+\tau\frac{d\tau}{d\alpha}\right)\\
 & = & \frac{l+\alpha\cdot \frac{l}{\alpha}(1+4\alpha^2 x)}{l^2+\tau^2}-\frac{2\alpha l}{(l^2+\tau^2)^2}\left(\frac{l^2}{\alpha}(1+4\alpha^2x)+\tau\left(\frac{\tau}{\alpha}+4\alpha y l\right)\right)\\
 & = & \frac 1{(l^2+\tau^2)^2}\left((2l+4\alpha^2lx)(l^2+\tau^2)-2l^3(1+4\alpha^2 x)-2\tau l(\tau +4\alpha^2 yl)\right)\\
 & = & \frac 1{(l^2+\tau^2)^2}\left(2l^3+2l\tau^2+4\alpha^2l^3x+4\alpha^2 \tau^2 lx-2l^3-8\alpha^2 l^3 x-2\tau^2l-8\tau\alpha^2yl^2\right)\\
  & = & \frac {4\alpha^2l}{(l^2+\tau^2)^2}\left(\tau^2 x- l^2 x-2\tau yl\right).
 \end{eqnarray*}
Applying the bounds on $|x|, |y|$ from theorem \ref{theorem:30} (3) gives:
$$\left|\frac{d}{dt}\left(\frac{\alpha \cdot l}{l^2+\tau^2}\right)\right|\leq \frac{2\alpha l}{(l^2+\tau^2)^2}\cdot \frac{1}{4\alpha^2 \sinh^2(R)}(\tau^2+l^2+2l\tau )\leq \frac{1}{2\alpha l  \sinh^2(R)}\left|\frac{l^2(l^2+2\tau l+\tau^2)}{(l^2+\tau^2)^2}\right|.$$
This is identical to the bound we attained for $\omega$ above. Thus, by an identical argument we may conclude that:
$$\left|\frac{m(0)\cdot h(0)}{h(0)^2+tw(0)^2}-\frac{2\pi \cdot l(2\pi)}{l(2\pi)^2+\tau(2\pi)^2}\right|\leq 92.$$
Again, considering figure \ref{figure:18}, we see that:
$$\frac{m(0)\cdot h(0)}{tw(0)^2+h(0)^2}=\frac{\text{Area}(\partial C)}{tw(0)^2+h(0)^2}=\frac{l(\gamma_\partial)\cdot \text{height}(\partial C)}{l(\gamma_\partial)^2}.$$
This completes the second inequality. $\square$

\newpage
\bibliographystyle{alpha}
\bibliography{citation}

\begin{thebibliography}{DDRW08}

\bibitem[Ago13]{Ago13}
Ian Agol.
\newblock The virtual {H}aken conjecture (with an appendix by {I}an {A}gol,
  {D}aniel {G}roves and {J}ason {M}anning).
\newblock {\em Documenta Mathematica}, 18:1045--1087, 2013.

\bibitem[B{\"o}r78]{Bor78}
Karoly B{\"o}r\"oczky.
\newblock Packing of spheres in spaces of constant curvature.
\newblock {\em Acta Mathematica Academia Scientarum Hungarica}, 32:243--261,
  1978.

\bibitem[Bro01a]{Bro01a}
Jeffrey Brock.
\newblock The {W}eil-{P}etersson metric and volumes of 3-dimensional hyperbolic
  convex cores.
\newblock {\em Journal of the American Mathematical Society}, 16(3):495--535,
  2001.

\bibitem[Bro01b]{Bro01b}
Jeffrey Brock.
\newblock {W}eil-{P}etersson translation distance and volumes of mapping tori.
\newblock {\em Communications in Analysis and Geometry}, 11(5):987--999, 2001.

\bibitem[Cap16]{Cap16}
Antonio~De Capua.
\newblock {\em Hyperbolic volume estimates via train tracks}.
\newblock PhD thesis, University of Oxford, 2016.

\bibitem[Cen22]{Cen22}
Mustafa Cengiz.
\newblock {H}eegard genus and complexity of fibered knots.
\newblock {\em Journal of Topology}, 15(4):2389--2425, 2022.

\bibitem[DDRW08]{DDRW08}
Patrick Dehornoy, Ivan Dynnikov, Dale Rolfson, and Bert Wiest.
\newblock {\em Ordering Braids}, volume 148 of {\em Mathematical Surverys and
  Monographs}.
\newblock American Mathematical Society, 2008.

\bibitem[Eli89]{Eli89}
Yakov Eliashberg.
\newblock Classification of overtwisted contact structures on 3-manifolds.
\newblock {\em Inventiones Mathematicae}, 98(3):623--637, 1989.

\bibitem[Etn04]{Etn04}
John Etnyre.
\newblock Lectures on open book decompositions and contact structures.
\newblock In David Ellwod, Peter Ozsv\'ath, Andr\'as Stipsicz, and Zolt\'an
  Szab\'o, editors, {\em {F}loer Homology, Guage Theory, and Low-Dimensional
  Topology}, pages 103--143. American Mathematical Society and Clay Mathematics
  Institute, 2004.

\bibitem[FKP08]{FKP08}
David Futer, Efstratia Kalfagianni, and Jessica Purcell.
\newblock {D}ehn filling, volume, and the {J}ones polynomial.
\newblock {\em Journal of Differential Geometry}, 78(3):429--464, 2008.

\bibitem[FPS22]{FPS22}
David Futer, Jessica Purcell, and Saul Schleimer.
\newblock Effective bilipschitz bounds on drilling and fillings.
\newblock {\em Geometry and Topology}, 26(3):1077--1188, 2022.

\bibitem[FS14]{FS14}
David Futer and Saul Schleimer.
\newblock Cusp geometry of fibered 3-manifolds.
\newblock {\em American Journal of Mathematics}, 136(2):309--356, 2014.

\bibitem[Gir02]{Gir02}
Emmanuel Giroux.
\newblock G\'eom\'etrie de contact: de la dimension trois vers les dimension
  sup\'erieures.
\newblock In {\em Proceedings of the International Congress of Mathematicians
  (Beijing 2002)}, volume~2, pages 405--414. Higher Education Press, Beijing,
  2002.

\bibitem[HK05]{HK05}
Craig Hodgson and Steven Kerckhoff.
\newblock Universal bounds for hyperbolic {D}ehn surgery.
\newblock {\em Annals of Mathematics}, 162(1):367--421, 2005.

\bibitem[HKM07]{HKM07}
Ko~Honda, William Kazez, and Gordana Mati\'c.
\newblock Right-veering diffeomorphisms of compact surfaces with boundary.
\newblock {\em Inventiones Mathematicae}, 169(2):427--449, 2007.

\bibitem[HKM08]{HKM08}
Ko~Honda, William Kazez, and Gordana Mati\'c.
\newblock Right-veering diffeomorphisms of compact surfaces with boundary ii.
\newblock {\em Geometry and Topology}, 12(4):2057--2094, 2008.

\bibitem[HM18]{HM18}
Matthew Hedden and Thomas Mark.
\newblock {F}loer homology and fractional {D}ehn twists.
\newblock {\em Advances in Mathematics}, 324:1--39, 2018.

\bibitem[IK17]{IK17}
Tetsuya Ito and Keiko Kawamuro.
\newblock Essential open book foliations and fractional {D}ehn twist
  coefficient.
\newblock {\em Geometriae Dedicata}, 187:17--67, 2017.

\bibitem[JS79]{JS79}
William Jaco and Peter Shalen.
\newblock {\em Seifert fibered spaces in 3-manifolds}, volume 220 of {\em
  Memoirs of the American Mathematical Society}.
\newblock American Mathematical Society, 1979.

\bibitem[Koj98]{Koj98}
Sodoyashi Kojima.
\newblock Deformations of hyperbolic 3-cone manifolds.
\newblock {\em Journal of Differential Geometry}, 49(3):469--516, 1998.

\bibitem[KR13]{KR13}
William Kazez and Rachel Roberts.
\newblock Fractional {D}ehn twists in knot theory and contact topology.
\newblock {\em Algebraic and Geometric Topology}, 13(6):3603--3637, 2013.

\bibitem[Lac19]{Lac19}
Marc Lackenby.
\newblock {D}ehn surgery from a hyperbolic perspective.
\newblock https://people.maths.ox.ac.uk/lackenby/, 2019.

\bibitem[Mag12]{Magid_Kleinian}
Aaron Magid.
\newblock Deformation spaces of {K}leinian surface groups are not locally
  connected.
\newblock {\em Geometry and Topology}, 16(3):1247--1320, 2012.

\bibitem[Mal05]{Mal05}
Andrei Malyutin.
\newblock {T}wist number of (closed) braids.
\newblock {\em St. Petersburg Mathematics Journal}, 16(5):791--813, 2005.

\bibitem[Mar22]{Mar22}
Bruno Martelli.
\newblock {\em An Introduction to {G}eometric {T}opology}.
\newblock arXiv.org, 2022.

\bibitem[Min10]{Min10}
Yair Minsky.
\newblock The classification of {K}leinian surface groups, i: models and
  bounds.
\newblock {\em Annals of Mathematics}, 171(1):1--107, 2010.

\bibitem[Mos68]{Mos68}
George Mostow.
\newblock Quasi-conformal mappings in n-space and the rigidity of hyperbolic
  space forms.
\newblock {\em Publications math\'ematiques de l'{I}.{H}.{\'E}.{S}.},
  34:53--104, 1968.

\bibitem[Pur08]{Purcell-Thesis}
Jessica Purcell.
\newblock Cusp shapes under cone deformation.
\newblock {\em Journal of Differential Geometry}, 80(3):453--500, 2008.

\bibitem[Thu79]{Thu79}
William Thurston.
\newblock {\em The geometry and topology of three-manifolds}.
\newblock Princeton Univ. Math. Dept. Notes, 1979.

\bibitem[Thu86]{Thu86}
William Thurston.
\newblock Hyperbolic structures on 3-manifolds, ii: Surfaces groups and
  3-manifolds which fiber over the circle.
\newblock {\em arXiv preprint}, 1998 (preprint first distributed in 1986).

\bibitem[Tol78]{Tol78}
Jeffrety Tollefson.
\newblock Involutions of {S}eifert fiber spaces.
\newblock {\em Pacific Journal of Mathematica}, 74(2):519--592, 1978.

\bibitem[Wal67]{Wal67}
Friedhelm Waldhausen.
\newblock Eine {K}lasse von 3-dimensionalen {M}annigfaltigkeiten ii.
\newblock {\em Inventiones Mathematicae}, 4:87--117, 1967.

\end{thebibliography}

\end{document}